\documentclass[a4paper,11pt]{article}
\usepackage{latexsym,amsfonts,amsmath}
\usepackage{graphics}

\def\inter{\mathop{\cap}}
\def\NN{\mathbb{N}}

\def\RR{\mathbb{R}}

\newcommand{\widebar}[1]{\overline{#1}}

\def\inter{\mathop{\cap}}
\def\liminf{\mathop{\underline{\lim}}}
\def\limsup{\mathop{\overline{\lim}}}
\def\ds{\displaystyle}

\def\nl{\mbox{} \newline }

\usepackage{natbib}

\begin{document}
\newtheorem{theorem}{Theorem}[section]
\newtheorem{proposition}[theorem]{Proposition}
\newtheorem{lemma}[theorem]{Lemma}
\newtheorem{corollary}[theorem]{Corollary}
\newtheorem{definition}[theorem]{Definition}
\newtheorem{remark}[theorem]{Remark}
\newtheorem{conjecture}[theorem]{Conjecture}
\newtheorem{assumption}[theorem]{Assumption}

\bibliographystyle{plain}

\title{The Vanishing Approach for the Average Continuous Control of Piecewise Deterministic Markov Processes
}
\author{ \mbox{ }
\\
O.L.V. Costa
\thanks{This author received financial support from CNPq (Brazilian National Research
Council), grant 304866/03-2 and FAPESP (Research Council of the
State of S\~ao Paulo), grant 03/06736-7.}
\\ \small Departamento de Engenharia de Telecomunica\c c\~oes e
Controle \\ \small Escola Polit\'ecnica da Universidade de S\~ao
Paulo \\ \small CEP: 05508 900-S\~ao Paulo, Brazil. \\ \small
e-mail: oswaldo@lac.usp.br
\\ \and
\\ F. Dufour
\thanks{Author to whom correspondence should be sent.}
\\
\small Universite Bordeaux I \\
\small IMB, Institut Math\'ematiques de Bordeaux \\
\small INRIA Bordeaux Sud Ouest, Team: CQFD \\
\small \small 351 cours de la Liberation \\
\small 33405 Talence Cedex, France \\ \small e-mail :
dufour@math.u-bordeaux1.fr }

\maketitle

\begin{abstract}
The main goal of this paper is to derive sufficient conditions for the existence of an optimal control strategy for the
long run average continuous control problem of piecewise deterministic Markov processes (PDMP's)
taking values in a general Borel space and with compact action space depending on the state variable. In order to do that we apply the so-called
vanishing discount approach (see \cite{hernandez96}, page 83) to obtain a solution to an average cost optimality inequality associated to the long
run average cost problem. Our main assumptions are written in terms of some integro-differential inequalities related to the so-called
expected growth condition, and geometric convergence of the post-jump location kernel associated to the PDMP.
\end{abstract}
\begin{tabbing}
\small \hspace*{\parindent}  \= {\bf Keywords:}
piecewise-deterministic Markov processes, continuous-time,
long-run average\\ cost, optimal control,
integro-differential optimality inequation, vanishing approach\\
\> {\bf AMS 2000 subject classification:} 60J25, 90C40, 93E20
\end{tabbing}

\newpage

\section{Introduction}
A general family of non-diffusion stochastic models suitable for formulating optimization problems in several areas of
operations research, namely piecewise-deterministic Markov processes (PDMP's), was introduced in \cite{davis84} and
\cite{davis93}. These processes are determined by three local characteristics; the flow $\phi$, the jump rate $\lambda$, and the
transition measure $Q$. Starting from $x$ the motion of the process follows the flow $\phi(x,t)$ until the first jump time
$T_1$ which occurs either spontaneously in a Poisson-like fashion with rate $\lambda$ or when the flow $\phi(x,t)$ hits the boundary
of the state-space. In either case the location of the process at the jump time $T_1$ is selected by the transition measure
$Q(\phi(x,T_1),.)$ and the motion restarts from this new point as before. A suitable choice of the state space and the local
characteristics $\phi$, $\lambda$, and $Q$ provide stochastic models covering a great number of problems of operations research
\cite{davis93}.

\bigskip

There exist two types of control for PDMP's: \textit{continuous control} and \textit{impulse control}. This terminology has been introduced by
M.H.A. Davis in \cite[page 134]{davis93} where continuous control is used to describe situations in which the control variable acts
at all times on the process through the characteristics $(\phi,\lambda,Q)$ by influencing the deterministic motion and the
probability of the jumps. On the other hand the terminology impulse control refers to a control that intervenes on the process by moving it to a new point of the state space at some times specifed by the controller.

\bigskip

In \cite{average} it was studied the long run average continuous control problem of PDMP's taking values in a general Borel space. At each point $x$ of the state space
a control variable is chosen from a compact action set $\mathbb{U}(x)$ and is applied on the jump parameter $\lambda$ and transition measure $Q$.
The goal was to minimize the long run average cost, which is composed of a running cost and a boundary cost (which is added each time the PDMP touches the boundary).
Both costs are assumed to be positive but not necessarily bounded. As far as the authors are aware of, this was the first time that this kind of problem was considered
in the literature. Indeed, results are available for the long run average cost problem but for impulse control see Costa
\cite{costa89}, Gatarek \cite{garatek93} and the book by M.H.A. Davis \cite{davis93} (see the references therein). On the other
hand, the continuous control problem has been studied only for discounted costs by A. Almudevar \cite{almudevar01}, M.H.A. Davis \cite{davis86,davis93}, M.A.H.
Dempster and J.J. Ye \cite{dempster92,dempster96}, Forwick, Sch{\"a}l, and Schmitz \cite{forwick04}, M. Sch\"al
\cite{schal98}, A.A. Yushkevich \cite{yushkevich87,yushkevich89}.

\bigskip

This paper deals with the vanishing approach for the long run average continuous control problem of a PDMP and can be seen as a continuation
of the results derived in \cite{average}.
By exploiting the special features of the PDMP's we trace a parallel with the general theory for discrete-time Markov Decision Processes (see, for instance, \cite{guo06,hernandez96}) rather than the
continuous-time case (see, for instance \cite{guo06a,Zhu08}).
The two main reasons for doing that is to use the powerful tools developed in the discrete-time framework
(see for example the references \cite{bertsekas78,dynkin79,hernandez96,hernandez99}) and to avoid working with the infinitesimal generator associated to a PDMP, which in most cases has its domain of definition
difficult to be characterized.
We develop further on the approach presented by the authors in \cite{average} which consists of using a connection between the continuous-time control problem of a PDMP and a discrete-time optimality equation
(see the introduction of section \ref{main3} for a detailed explanation of this method).
In particular, we derive sufficient conditions under which a
boundedness condition (with the lower bound being a function rather than a constant as supposed in \cite{average}) on the value functions for the discounted problems is satisfied.
The main assumptions for this are based on some integro-differential inequalities related to the so-called expected growth condition (see Assumption \ref{A1}), and geometric
convergence of the post-jump location kernel associated to the PDMP (see Assumption \ref{A2}).
As a consequence, we obtain a result of existence of an optimal ordinary control strategy for the long run average control problem of a PDMP having the important property of being in a feedback form.

\bigskip

The paper is organized in the following way. In section \ref{Intro} we introduce some notation, basic assumptions, and the problem formulation.
In section \ref{AuxRes} we introduce several assumptions related to the continuity of the parameters, the expected growth condition and geometric convergence of the
post-jump location of the PDMP. In the sequence we provide several key auxiliary results for obtaining a bound for the discounted problems, and some extensions of the results
presented in \cite{average} to the case in which the functions under consideration are not necessarily positive but just bounded by a test function $g$.
The main results are presented in section \ref{main3}, which provides sufficient conditions for the existence of an optimal control strategy for the
long run average continuous control problem of a PDMP and obtain a solution to an average cost optimality inequality associated to the long
run average cost problem.

\section{Notation, basic assumptions, and problem formulation}\label{Intro}
\subsection{Presentation of the control problem}
\label{pre}
In this section we present some standard notation and some basic definitions related to the motion of a PDMP $\{X(t)\}$,
and the control problems we will consider throughout the paper.
For further details and properties the reader is referred to \cite{davis93}. The following notation will be used in this paper:
$\NN$ denotes the set of natural numbers, $\RR$ the set of real numbers, $\RR_+$ the set of positive real numbers and $\RR^d$ the $d$-dimensional euclidian space. We write
$\eta$ as the Lebesgue measure on $\RR$. For $X$ a metric space $\mathcal{B}(X)$ represents the $\sigma$-algebra generated by the open sets of $X$. $\mathcal{M}(X)$ (respectively, $\mathcal{P}(X)$)
denotes the set of all finite (respectively probability) measures on $(X,\mathcal{B}(X))$. Let $X$ and $Y$ be metric spaces.
The set of all Borel measurable (respectively bounded) functions from $X$ into $Y$ is denoted by $\mathbb{M}(X;Y)$ (respectively $\mathbb{B}(X;Y)$).
Moreover, for notational simplicity $\mathbb{M}(X)$ (respectively $\mathbb{B}(X)$, $\mathbb{M}(X)^{+}$, $\mathbb{B}(X)^{+}$) denotes $\mathbb{M}(X;\RR)$
(respectively $\mathbb{B}(X;\RR)$, $\mathbb{M}(X;\RR_{+})$, $\mathbb{B}(X;\RR_{+})$).
For $g\in \mathbb{M}(X)$ with $g(x)>0$ for all $x\in X$, $\mathbb{B}_{g}(X)$ is the set of functions  $v\in \mathbb{M}(X)$ such that $\ds ||v(x)||_{g}=\sup_{x\in X} \frac{|v(x)|}{g(x)}< +\infty$.
$\mathbb{C}(X)$ denotes the set of continuous functions from $X$ to $\RR$.
For $h\in \mathbb{M}(E)$, $h^+$ (respectively $h^{-}$) denotes the positive (respectively, negtive) part of $h$.

\noindent
Let $E$ be an open subset of $\RR^n$, $\partial E$ its boundary, and $\widebar{E}$ its closure. A controlled PDMP is determined by its local
characteristics $(\phi,\lambda,Q)$, as presented in the sequel.
The flow $\phi(x,t)$ is a function $\phi: \: \mathbb{R}^{n}\times \RR_{+} \longrightarrow \mathbb{R}^{n}$ continuous in $(x,t)$ and such that $\phi(x,t+s) = \phi(\phi(x,t),s).$
For each $x\in E$ the time the flow takes to reach the boundary starting from $x$ is defined as
$t_{*}(x)\doteq \inf \{t>0:\phi(x,t)\in \partial E \}$. For $x\in E$ such that $t_{*}(x)=\infty$ (that is, the flow starting from $x$ never touches the boundary), we set
$\phi(x,t_{*}(x))=\Delta$, where $\Delta$ is a fixed point in $\partial E$.
 We define the following space of functions absolutely continuous along the flow with limit towards the boundary:
\begin{align*}
\mathbb{M}^{ac}(E) & =  \bigl\{ g\in\mathbb{M}(E) \: : \: g(\phi(x,t)): [0,t_{*}(x)) \mapsto \RR \text{ is absolutely continuous for each } x\in E\\
& \text{ and whenever } t_{*}(x)< \infty \text{ the limit }\lim_{t\rightarrow t_{*}(x)} g(\phi(x,t)) \text{ exists} \bigr\}.
\end{align*}
For $g\in \mathbb{M}^{ac}(E)$ and $z\in \partial E$ for which there exists $x\in E$ such that $z=\phi(x,t_{*}(x))$ where $t_{*}(x)< \infty$ we define $\ds g(z) = \lim_{t\rightarrow t_{*}(x)} g(\phi(x,t))$
(note that the limit exists by assumption). As shown in Lemma 2 in \cite{ECC07}, for $g\in\mathbb{M}^{ac}(E)$ there exists a function $\mathcal{X}g \in \mathbb{M}(E)$ such that for all $x\in E$ and $t\in [0,t_{*}(x))$
$g(\phi(x,t))-g(x)  =  \int_{0}^{t} \mathcal{X}g(\phi(x,s)) ds $.

\bigskip

\noindent
The local characteristics $\lambda$ and $Q$ depend on a control action $u\in \mathbb{U}$ where $\mathbb{U}$ is a compact metric space
(there is no loss of generality in assuming this property for $\mathbb{U}$, see Remark 2.8 in \cite{average}), in the following way:
$\lambda \in \mathbb{M}(\widebar{E}\times\mathbb{U})^{+}$ and 
$Q$ is a stochastic kernel on $E$ given $\widebar{E}\times \mathbb{U}$.
For each $x\in \widebar{E}$ we define the subsets $\mathbb{U}(x)$ of $\mathbb{U}$ as the set of feasible control actions that can be taken when the state process is in $x\in \widebar{E}$,
that is, the control action that will be applied to $\lambda$ and $Q$ must belong to $\mathbb{U}(x)$.
The following assumptions, based on the standard theory of Markov decision processes (see for example \cite{hernandez96}), will be made throughout the paper:
\begin{assumption}
\label{Hyp1a} For all $x\in \widebar{E}$, $\mathbb{U}(x)$ is a compact subspace of $\mathbb{U}$.
\end{assumption}
\begin{assumption}
\label{Hyp2a} The set $K=\left\{(x,a): x\in \widebar{E}, a \in \mathbb{U}(x) \right\}$ is a Borel subset of $\widebar{E}\times \mathbb{U}$.
\end{assumption}

\noindent We present next the definition of an admissible control strategy and the associated motion of the controlled process.
A control policy $U$ is a pair of functions $(u,u_{\partial}) \in \mathbb{M}(\NN \times E\times \RR_{+};\mathbb{U}) \times \mathbb{M}(\NN \times E;\mathbb{U})$
satisfying $u(n,x,t)\in \mathbb{U}(\phi(x,t))$, and $u_{\partial}(n,x)\in \mathbb{U}(\phi(x,t_{*}(x)))$ for all $(n,x,t)\in \NN \times E\times \RR_{+}$.
The class of admissible control strategies will be denoted by $\mathcal{U}$.
Consider the state space $\widehat{E}=E\times E \times \RR_{+} \times \NN$.
For a control policy $U=(u,u_{\partial})$ let us introduce the following parameters for $\hat{x}=(x,z,s,n)\in \widehat{E}$:
the flow $\widehat{\phi}(\hat{x},t) = (\phi(x,t), z , s+t , n)$,
the jump rate $\widehat{\lambda}^{U}(\hat{x})=\lambda(x,u(n,z,s))$, and
the transition measure
\begin{eqnarray*}
\widehat{Q}^{U}(\hat{x},A\times B\times \{0\}\times \{n+1\}) =
\begin{cases}
Q(x,u(n,z,s)); A\inter B) & \text{ if } x\in E,\\
Q(x,u_{\partial}(n,z);A\inter B) & \text{ if } x\in \partial E,
\end{cases}
\end{eqnarray*}
for $A$ and $B$ in $\mathcal{B}(E)$. From \cite[section 25]{davis93}, it can be shown that for any control strategy $U=(u,u_{\partial})\in \mathcal{U}$
there exists a filtered probability space $(\Omega,\mathcal{F},\{ \mathcal{F}_{t} \}, \{ P^{U}_{\hat{x}} \}_{\hat{x}\in \widehat{E}})$
such that the piecewise deterministic Markov process $\{\widehat{X}^{U}(t)\}$ with local characteristics $(\widehat{\phi},\widehat{\lambda}^{U},\widehat{Q}^{U})$ may be constructed as follows.
For notational simplicity the probability $P^{U}_{\hat{x}_{0}}$ will be denoted by $P^{U}_{(x,k)}$ for $\hat{x}_{0}=(x,x,0,k)\in \widehat{E}$.
Take a random variable $T_1$ such that
\begin{equation*}
P^{U}_{(x,k)}(T_1>t) \doteq
\begin{cases}
e^{-\Lambda^{U}(x,k,t)} & \text{for } t<t_{*}(x)\\
0 & \text{for } t\geq t_{*}(x)
\end{cases}
\end{equation*}
where for $x\in E$ and $t\in [0,t_*(x)[$, $\Lambda^{U}(x,k,t) \doteq \int_0^t\lambda(\phi(x,s),u(k,x,s))ds.$
If $T_1$ is equal to infinity, then for $t\in \RR_+$, $\widehat{X}^{U}(t)= \bigl(\phi(x,t),x,t,k\bigr)$.
Otherwise select independently an $\widehat{E}$-valued random variable (labelled $\widehat{X}^{U}_{1}$) having distribution
\begin{align*}
P^{U}_{(x,k)}(\widehat{X}^{U}_{1} \in A \times B \times \{0\}
\times \{k+1\} |\sigma\{T_1\}) =
\begin{cases}
Q(\phi(x,T_{1}),u(k,x,T_{1})); A\inter B) & \text{ if } \phi(x,T_{1})\in E, \\
Q(\phi(x,T_{1}),u_{\partial}(k,x); A\inter B) & \text{ if }
\phi(x,T_{1}) \in \partial E.
\end{cases}
\end{align*}
The trajectory of $\{\widehat{X}^{U}(t)\}$ starting from $(x,x,0,k)$, for $t\leq T_1$ , is given by
\begin{equation*}
\widehat{X}^{U}(t) \doteq
\begin{cases}
\bigl(\phi(x,t),x,t,k\bigr) &\text{for } t<T_1, \\
\widehat{X}^{U}_1 &\text{for } t=T_1.
\end{cases}
\end{equation*}
Starting from $\widehat{X}^{U}(T_1)=\widehat{X}^{U}_1$, we now select the next inter-jump time $T_2-T_1$ and post-jump location $\widehat{X}^{U}(T_2)=\widehat{X}^{U}_2$ in a similar way.
Let us define the components of the PDMP $\{\widehat{X}^{U}(t)\}$ by $\widehat{X}^{U}(t)=\bigl(X(t),Z(t),\tau(t),N(t)\bigr)$.
From the previous construction, it is easy to see that $X(t)$ corresponds to the trajectory of the
system, $Z(t)$ is the value of $X(t)$ at the last jump time before $t$, $\tau(t)$ is time elapsed between the last jump and time $t$,
and $N(t)$ is the number of jumps of the process  $\{X(t)\}$ at time $t$. As in Davis \cite{davis93}, we consider the following assumption to avoid any accumulation point of the jump times:
\begin{assumption}
\label{Hypjump} For any $x\in E$, $U=(u,u_{\partial})\in
\mathcal{U}$, and $t\geq 0$, we have $\ds E^{U}_{(x,0)}\Biggl[
\sum_{i=1}^{\infty} I_{\{T_{i}\leq t\}} \Biggr] < \infty$.
\end{assumption}

\noindent The costs of our control problem will contain two terms, a running cost $f$ and a boundary cost $r$, satisfying the following properties:
\begin{assumption}
\label{Hyp6a} $f\in \mathbb{M}(\widebar{E}\times\mathbb{U})^{+}$, and $r\in \mathbb{M}(\partial E\times\mathbb{U})^{+}$.
\end{assumption}
Define for $\alpha \geq 0$, $t\in \mathbb{R}_+$, and $U\in \mathcal{U}$,
\begin{align*}
\mathbf{J}^\alpha(U,t) = \int_{0}^{t} e^{-\alpha s}f\bigl(X(s), & u(N(s),Z(s),\tau(s)) \bigr) ds  + \int_{0}^{t}e^{-\alpha s} r\bigl(X(s-),u_{\partial}(N(s-),Z(s-)) \bigr)
dp^{*}(s),
\end{align*}
where $\ds p^{*}(t) = \sum_{i=1}^{\infty} I_{\{T_{i}\leq t\}} I_{\{ X(T_{i}-) \in \partial E\}}$ counts the number of times the
process hits the boundary up to time $t$ and, for notational simplicity, set $\mathbf{J}(U,t)=\mathbf{J}^0(U,t)$.
The long-run average cost we want to minimize over $\mathcal{U}$ is given by:
$\mathcal{A}(U,x) = \limsup_{t\rightarrow +\infty} \frac{1}{t} E^{U}_{(x,0)} [\mathbf{J}(U,t)]$
and we set $\mathcal{J}_{\mathcal{A}}(x) = \inf_{U\in\mathcal{U}}\mathcal{A}(U,x)$.
For the $\alpha$ discounted case, with $\alpha > 0$, the cost we want to minimize is given by:
$\mathcal{D}^{\alpha}(U,x) = E^{U}_{(x,0)}[\mathbf{J}^\alpha(U,\infty)]$ and we set
$\mathcal{J}_{\mathcal{D}}^{\alpha}(x) = \inf_{U\in\mathcal{U}}\mathcal{D}^{\alpha}(U,x)$.
We need the following assumption, to avoid infinite costs for the discounted case.
\begin{assumption}
\label{Hypdis} For all $\alpha>0$ and all $x\in E$,
$\mathcal{J}_{\mathcal{D}}^{\alpha}(x)<\infty$.
\end{assumption}

\subsection{Discrete-time relaxed and ordinary controls}
\label{conrex1}
We present in this sub-section the set of discrete-time relaxed and ordinary controls.
\nl
Consider $\mathbb{C}(\mathbb{U})$ equipped with the topology of uniform convergence and $\mathcal{M}(\mathbb{U})$ equipped with the
weak$^{*}$ topology $\sigma(\mathcal{M}(\mathbb{U}),\mathbb{C}(\mathbb{U}))$.
For $x\in E$, define $\mathcal{P}_{x}\bigl(\mathbb{U}\bigr)$ as the set of measures $\mu\in\mathcal{P}(\mathbb{U})$ satisfying $\mu(\mathbb{U}(\phi(x,t_{*}(x))))=1$.
$\mathcal{P}(\mathbb{U})$ and $\mathcal{P}_{x}(\mathbb{U})$ for $x\in E$ are subsets of $\mathcal{M}(\mathbb{U})$ and are equipped with the relative topology.

Let $\mathcal{V}^{r}$ (respectively $\mathcal{V}^{r}(x)$ for $x\in E$) be the set of all $\eta$-measurable functions $\mu$ defined on $\RR_{+}$ with value in $\mathcal{P}(\mathbb{U})$ such that
$\mu(t,\mathbb{U})=1$ $\eta$-a.e. (respectively $\mu(t,\mathbb{U}(\phi(x,t)))=1$ $\eta$-a.e.).
It can be shown (see sub-section 3.1 in \cite{average}) that $\mathcal{V}^{r}(x)$ is a compact set of the metric space $\mathcal{V}^{r}$: a sequence $\bigl(\mu_{n}\bigr)_{n\in \NN}$ in $\mathcal{V}^{r}(x)$ converges to $\mu$ if and only if for all $g\in L^{1}(\RR_{+};\mathbb{C}(\mathbb{U}))$
\begin{eqnarray*}
\lim_{n\rightarrow \infty} \int_{\RR_{+}} \int_{\mathbb{U}(\phi(x,t))}  g(t,u) \mu_{n}(t,du) dt = \int_{\RR_{+}} \int_{\mathbb{U}(\phi(x,t))}  g(t,u) \mu(t,du) dt.
\end{eqnarray*}
The sets of relaxed controls can be defined as follows:
$\mathbb{V}^{r}(x)  = \mathcal{V}^r(x) \times \mathcal{P}_{x}\bigl(\mathbb{U}\bigr)$, for $x\in E$ 
and $\mathbb{V}^{r}  = \mathcal{V}^r \times \mathcal{P}\bigl(\mathbb{U}\bigr)$.
The set of ordinary controls, denoted by $\mathbb{V}$ (respectively $\mathbb{V}(x)$ for $x\in E$), is defined as above
except that it is composed of deterministic functions instead of probability measures. More specifically we have
$\mathcal{V}(x)  =  \bigl\{ \nu \in \mathbb{M}(\RR_{+}, \mathbb{U}) : (\forall t\in \RR_{+}), \nu(t) \in \mathbb{U}(\phi(x,t)) \bigr\}$,
$\mathbb{V}(x)  =  \mathcal{V}(x) \times \mathbb{U}(\phi(x,t_{*}(x))),$
$\mathbb{V}  =  \mathbb{M}(\RR_{+}, \mathbb{U}) \times
\mathbb{U}$.
Consequently, the set of ordinary controls is a subset of the set of relaxed controls $\mathbb{V}^{r}$ (respectively
$\mathbb{V}^{r}(x)$ for $x\in E$) by identifying any control action $u\in \mathbb{U}$ with the Dirac measure concentrated on
$u$. Thus we can write that $\mathbb{V}\subset \mathbb{V}^{r}$ (respectively $\mathbb{V}(x)\subset \mathbb{V}^{r}(x)$ for $x\in
E$) and from now on we will consider that $\mathbb{V}$ (respectively $\mathbb{V}(x)$ for $x\in E$) will be endowed with
the topology generated by $\mathbb{V}^{r}$.
The necessity to introduce the class of relaxed control $\mathbb{V}^{r}$ is justified by the fact that in general there
does not exist a topology for which $\mathbb{V}$ and $\mathbb{V}(x)$ are compact sets.

\noindent
As in \cite{hernandez96}, page 14, we need that the set of feasible state/relaxed-control pairs is a measurable subset of $\mathcal{B}(E)\times \mathcal{B}(\mathbb{V}^{r})$, that is, we need the following assumption.
\begin{assumption}
\label{Mesurability}
$\mathcal{K} \doteq \bigl\{ (x,\Theta) : \Theta \in \mathbb{V}^{r}(x), x\in E \bigr\} \in \mathcal{B}(E)\times \mathcal{B}(\mathbb{V}^{r}).$
\end{assumption}
A sufficient condition is presented in \cite[Proposition 3.3]{average} to ensure that Assumption \ref{Mesurability} holds.

\subsection{Discrete-time operators and measurability properties}
In this sub-section we present some important operators associated to the optimality equation of the discrete-time problem.
We consider the following notation $\ds w(x,\mu)  \doteq  \int_{\mathbb{U}} w(x,u)  \mu (du)$ and
$\ds Qh(x,\mu) \doteq \int_{\mathbb{U}} \int_{E} h(z) Q(x,u;dz) \mu (du)$, and
$\ds \lambda Qh(x,\mu) \doteq  \int_{\mathbb{U}} \lambda(x,u) \int_{E} h(z) Q(x,u;dz)  \mu (du)$
for $x\in \widebar{E}$, $\mu\in\mathcal{P}\bigl(\mathbb{U}\bigr)$, $h\in \mathbb{M}(E)^{+}$ and $w\in \mathbb{M}(\widebar{E}\times \mathbb{U})^{+}$.

The following operators will be associated to the optimality equations of the discrete-time problems that will be
presented in the next sections. For $\Theta=\bigl(\mu,\mu_{\partial}\bigr)\in \mathbb{V}^{r}$, $(x,A)\in E\times \mathcal{B}(E)$, $\alpha \in \RR$, according to Lemma 2 in \cite[Appendix 5]{dynkin79}
define
\begin{eqnarray}
\Lambda^{\mu}(x,t)  & \doteq & \int_{0}^{t} \lambda(\phi(x,s),\mu (s)) ds \nonumber \\
\label{DefGr}
G_{\alpha}(x,\Theta;A) & \doteq & \int_0^{t_{*}(x)}e^{-\alpha s - \Lambda^{\mu}(x,s)}\lambda QI_{A}(\phi(x,s),\mu(s)) ds\nonumber \\
& &+ e^{-\alpha t_{*}(x) -\Lambda^{\mu}(x,t_{*}(x))} Q(\phi(x,t_{*}(x)),\mu_{\partial};A).
\end{eqnarray}
For $h\in \mathbb{M}(E)^{+}$, we define $G_{\alpha}h(x,\Theta)\doteq\ds \int_{E} h(y) G_{\alpha}(x,\Theta;dy)$.
For $x\in E$, $\Theta=\bigl(\mu,\mu_{\partial}\bigr)\in \mathbb{V}^{r}$, $v\in \mathbb{M}(E\times \mathbb{U})^{+}$, $w\in \mathbb{M}(\partial E\times \mathbb{U})^{+}$,
$\alpha \in \RR$, introduce
\begin{eqnarray}
\label{DefLr}
L_{\alpha}v(x,\Theta) & \doteq  & \int_0^{t_{*}(x)}e^{-\alpha s-\Lambda^{\mu}(x,s)} v(\phi(x,s),\mu(s)) ds, \\
\label{DefHr} H_{\alpha}w(x,\Theta) & \doteq & e^{-\alpha t_{*}(x)-\Lambda^{\mu}(x,t_{*}(x))} w(\phi(x,t_{*}(x)),\mu_{\partial}).
\end{eqnarray}

For $h\in \mathbb{M}(E)$ (respectively, $v\in \mathbb{M}(E\times \mathbb{U})$),  $G_{\alpha}h(x,\Theta)=G_{\alpha}h^{+}(x,\Theta)-G_{\alpha}h^{-}(x,\Theta)$
(respectively, $L_{\alpha}v(x,\Theta)=L_{\alpha}v^{+}(x,\Theta)-L_{\alpha}v^{-}(x,\Theta)$) provided the difference has a meaning.
It will be useful in the sequel to define the function $\mathcal{L}_{\alpha}(x,\Theta)$  as follows:
$\mathcal{L}_{\alpha}(x,\Theta)  \doteq  L_{\alpha}I_{E\times \mathbb{U}}(x,\Theta)$. In particular for $\alpha =0$ we write for simplicity $G_0=G$, $L_0=L$, $H_0=H$, $\mathcal{L}_0=\mathcal{L}$.
Measurability properties of the operators $G_{\alpha}$, $L_{\alpha}$, and $H_{\alpha}$ are shown in \cite[Proposition 3.4]{average}.
\bigskip

\noindent We present now the definitions of the one-stage optimization operators.
\begin{definition} Let $\alpha\in \RR_{+}$, $\rho\in \RR$, and $h\in \mathbb{M}(E)$.
Assume that for any $x\in E$ and $\Upsilon\in \mathbb{V}(x)$, $-\rho \mathcal{L}_{\alpha}(x,\Upsilon) +L_{\alpha}f(x,\Upsilon)+H_{\alpha}r(x,\Upsilon)+G_{\alpha} h(x,\Upsilon) $ is well defined.
The (ordinary) one-stage optimization operator is defined by
\begin{equation*}
\mathcal{T}_{\alpha}(\rho,h)(x)  =  \inf_{\Upsilon\in \mathbb{V}(x)} \Bigl\{-\rho \mathcal{L}_{\alpha}(x,\Upsilon)
+L_{\alpha}f(x,\Upsilon)+H_{\alpha}r(x,\Upsilon)+G_{\alpha} h(x,\Upsilon)  \Bigr\}.
\end{equation*}
Assume that for any $x\in E$ and $\Theta\in \mathbb{V}^{r}(x)$, $-\rho \mathcal{L}_{\alpha}(x,\Theta) +L_{\alpha}f(x,\Theta)+H_{\alpha}r(x,\Theta)+G_{\alpha} h(x,\Theta)$ is well defined.
The relaxed one-stage optimization operator is defined by
\begin{equation*}
\mathcal{R}_{\alpha}(\rho,h)(x)  = \inf_{\Theta\in \mathbb{V}^{r}(x)} \Bigl\{-\rho \mathcal{L}_{\alpha}(x,\Theta)
+L_{\alpha}f(x,\Theta)+H_{\alpha}r(x,\Theta)+G_{\alpha} h(x,\Theta) \Bigr\}.
\end{equation*}
\end{definition}

\noindent In particular for $\alpha =0$ we write for simplicity $\mathcal{T}_{0}=\mathcal{T}$, and $\mathcal{R}_{0}=\mathcal{R}$.

\bigskip

\noindent The sets of measurable selectors associated to $\bigl(\mathbb{U}(x)\bigr)_{x\in E}$, $\bigl(\mathbb{V}(x)\bigr)_{x\in E}$, $\bigl(\mathbb{V}^{r}(x)\bigr)_{x\in E}$ are defined by
$\mathcal{S}_{\mathbb{U}}  =  \Bigl\{  u \in \mathbb{M}(\widebar{E}, \mathbb{U}) : (\forall x\in \widebar{E}), u(x) \in \mathbb{U}(x)\Bigr\}$,
$\mathcal{S}_{\mathbb{V}}  =  \Bigl\{ (\nu,\nu_{\partial})\in \mathbb{M}(E, \mathbb{V}) : (\forall x\in E), \bigl( \nu(x), \: \nu_{\partial}(x) \bigr)\in \mathbb{V}(x)\Bigr\}$,
$\mathcal{S}_{\mathbb{V}^{r}} = \Bigl\{ (\mu,\mu_{\partial})\in \mathbb{M}(E, \mathbb{V}^{r}) : (\forall x\in E), \bigl( \mu(x), \: \mu_{\partial}(x) \bigr)\in \mathbb{V}^{r}(x)\Bigr\}$.

\bigskip

\noindent
\noindent For $\alpha\in \RR_{+}$, $\rho\in \RR$, and $v\in \mathbb{M}(E)$, the one-stage optimization problem associated to the operator $\mathcal{T}_{\alpha}(\rho,v)$, respectively
$\mathcal{R}_{\alpha}(\rho,v)$, consists of finding a measurable selector $\Upsilon\in \mathcal{S}_{\mathbb{V}}$, respectively $\Theta\in \mathcal{S}_{\mathbb{V}^{r}}$ such that for all $x\in E$,
$\mathcal{T}_{\alpha}(\rho,v)(x)  =  -\rho \mathcal{L}_{\alpha}(x,\Upsilon) + L_{\alpha}f(x,\Upsilon)+H_{\alpha}r(x,\Upsilon)+G_{\alpha} v(x,\Upsilon)$ and respectively
$\mathcal{R}_{\alpha}(\rho,v)(x)  =  -\rho \mathcal{L}_{\alpha}(x,\Theta) + L_{\alpha}f(x,\Theta)+H_{\alpha}r(x,\Theta)+G_{\alpha} v(x,\Theta)$.

\bigskip

\noindent Finally we conclude this section by recalling (see Propositions 3.8 and 3.10 in \cite{average}) that there exist two natural mappings from $\mathcal{S}_{\mathbb{U}}$ to $\mathcal{S}_{\mathbb{V}}$ and from $\mathcal{S}_{\mathbb{U}}$ to $\mathcal{U}$.
\begin{definition}
\label{mapu}
For $u \in\mathcal{S}_{\mathbb{U}}$, define the measurable mapping $u_{\phi}$ of the space $E$ into $\mathbb{V}$ by
\nl
$u_{\phi}$ $:$ $x$ $\rightarrow$ $\bigl(u(\phi(x,.)),u(\phi(x,t_{*}(x)))\bigr)$.
\end{definition}

\begin{definition}
\label{mapU}
For $u \in\mathcal{S}_{\mathbb{U}}$, define the measurable mapping $U_{u_{\phi}}$ of the space $\NN\times E\times \RR_{+}$ into $\mathbb{U}\times \mathbb{U}$
by $U_{u_{\phi}}$ $:$ $(n,x,t)$ $\rightarrow$ $\bigl(u(\phi(x,t)),u(\phi(x,t_{*}(x)))\bigr)$ of the space $\NN\times E\times \RR_{+}$ into $\mathbb{U}\times \mathbb{U}$.
\end{definition}

\begin{remark}
\label{feedbak} The measurable selectors of the kind $u_{\phi}$ as in Definition \ref{mapu} are called ordinary feedback measurable selectors in the class $\mathcal{S}_{\mathbb{V}}\subset \mathcal{S}_{\mathbb{V}^{r}}$
and the control strategies of the kind $U_{u_{\phi}}$ as in definition \ref{mapU} are called ordinary feedback control strategies in the class $\mathcal{U}$.
\end{remark}

\section{Assumptions and auxiliary results}\label{AuxRes}
The purpose of this section is to introduce several assumptions (see sub-section \ref{AssDef}) and to derive preliminary results that will ensure the existence of an optimal control for the long run average cost.
More specifically, the two main results of sub-section \ref{Value} consist, roughly speaking, of providing a bound for $\mathcal{J}_D^\alpha(x)$ in terms of $\alpha$ (see Corollary \ref{JD})
and of proving that the mapping defined by $\mathcal{J}_D^\alpha(\cdot)$-$\mathcal{J}_D^\alpha(y)$ for $y$ fixed in $E$ belongs to $\mathbb{B}_{g}(E)$ (see  Theorem \ref{theo2}).
The results of sub-section \ref{auxil} are extensions of those presented in \cite{average} to the case in which the functions under consideration are not necessarily positive (as it was supposed in
\cite{average}) but instead belong to $\mathbb{B}_{g}(E)$. It must be pointed out that these generalizations are not straightforward and are crucial for obtaining the existence of an optimal ordinary feedback
control strategy for the long run average-cost problem of a PDMP.
In particular, Theorem \ref{theo3main} states that for any function $h\in \mathbb{B}_{g}(E)$, the one-stage optimization operators $\mathcal{R}_{\alpha}(\rho,h)(x)$ and $\mathcal{T}_{\alpha}(\rho,h)(x)$ are equal and
that there exists an ordinary feedback measurable selector for the one-stage optimization problems associated to these operators.

\subsection{Assumptions and definitions}
\label{AssDef}
The next assumption is somehow related to the so-called expected growth condition (see, for instance, Assumption 3.1 in \cite{guo06} for the discrete-time case, or
Assumption A in \cite{guo06a} for the continuous-time case) used, among other things, to guarantee uniform boundedness of $\mathcal{J}_D^\alpha(x)$ with respect to $\alpha$.

\begin{assumption}
\label{A1}
Suppose that there exist $b\geq 0$, $c> 0$, $\delta>0$, $M\geq 0$ and $g\in\mathbb{M}^{ac}(E)$, $g\geq 1$
$\overline{r} \in\mathbb{M}(\partial E)$, $\overline{r}(z)\geq 0$,
satisfying for all $x\in E$
\begin{eqnarray}
& \ds \sup_{a\in \mathbb{U}(x)}\Bigl\{\mathcal{X}g(x)+c g(x)-\lambda(x,a)\left[g(x)-Qg(x,a)\right]\Bigr\} \leq b, &
\label{Cu1} \\
& \ds \sup_{a\in \mathbb{U}(x)}\Bigl\{f(x,a)\Bigr\}\leq M g(x), &
\label{Cu3}
\end{eqnarray}
and for all $x \in E$ with $t_{*}(x)<\infty$
\begin{eqnarray}
& \ds \sup_{a\in \mathbb{U}(\phi(x,t_{*}(x)))}\{\overline{r}(\phi(x,t_{*}(x)))+Qg(\phi(x,t_{*}(x)),a)\} \leq g(\phi(x,t_{*}(x))), &
\label{Cu2} \\
& \ds \sup_{a\in \mathbb{U}(\phi(x,t_{*}(x)))}\Bigl\{r(\phi(x,t_{*}(x)),a)\Bigr\}\leq \frac{M}{c+\delta}\overline{r} (\phi(x,t_{*}(x))).&
\label{Cu3a}
\end{eqnarray}
\end{assumption}
Assumptions \ref{Hyp3a}, \ref{Hyp6bis} and \ref{Hyp5a}, presented in the sequel, are needed to guarantee some convergence and semi-continuity properties of the one-stage optimization operators (see sub-section
\ref{auxil}), and the existence of a measurable selector.
\begin{assumption}
\label{Hyp3a} For each $x\in E$, the restriction of $\lambda(x,.)$ to $\mathbb{U}(x)$ is continuous, for $t\in [0,t_{*}(x))$, $\ds \int_{0}^{t} \sup_{a\in \mathbb{U}(\phi(x,s))} \lambda(\phi(x,s),a) \: ds < \infty$ and if $t_{*}(x)< \infty$ then $\ds \int_{0}^{t_{*}(x)} \sup_{a\in \mathbb{U}(\phi(x,s))} \lambda(\phi(x,s),a) \: ds < \infty$.
\end{assumption}

\begin{assumption}
\label{Hyp6bis} There exists a sequence of measurable functions
$(f_j)_{j\in \NN}$ in $\mathbb{M}(\widebar{E}\times \mathbb{U})^{+}$ such that for all $y\in \widebar{E}$, $f_j(y,.)\uparrow f(y,.)$ as $j\rightarrow \infty$ and the restriction of $f_j(y,.)$ to $\mathbb{U}(y)$
is continuous.
There exists a sequence of measurable functions
$(r_j)_{j\in \NN}$ in $\mathbb{M}(\partial E\times \mathbb{U})^{+}$ such that for all $z\in \partial E$,
$r_j(z,.)\uparrow r(z,.)$ as $j\rightarrow \infty$ and the restriction of $r_j(z,.)$ to $\mathbb{U}(z)$ is continuous.\end{assumption}
\begin{assumption}
\label{Hyp5a} For all $x \in \widebar{E}$ and $h\in \mathbb{B}(E)$, the restriction of $Qh(x,.)$ to $\mathbb{U}(x)$ is continuous.
\end{assumption}
We make the following definition:
\begin{definition}
\label{defmes}
Consider $w\in\mathbb{M}(E)$ and $h\in\mathbb{B}_{g}(E)$. We define:
\begin{enumerate}
\item [D1)] $\widehat{u}(w,h)\in \mathcal{S}_{\mathbb{U}}$ as the measurable selector satisfying
\begin{align*}
\inf_{a\in  \mathbb{U}(x)} \{f(x,a) -\lambda(x, & a)  \Bigl[ w(x)-Q h(x,a) \Bigr]\} \\
& =  f(x,\widehat{u}(w,h)(x))-\lambda(x,\widehat{u}(w,h)(x)) \Bigl[w(x)-Qh(x,\widehat{u}(w,h)(x)) \Bigr],
\end{align*}
\begin{eqnarray*}
\inf_{a\in \mathbb{U}(z)}\{r(z,a)+Qh(z,a)\} & = & r(z,\widehat{u}(w,h)(z))+Qh(z,\widehat{u}(w,h)(z)).
\end{eqnarray*}
\item [D2)] $\widehat{u}_{\phi}(w,h)\in \mathcal{S}_{\mathbb{V}}$ as the measurable selector derived from $\widehat{u}(w,h)$ through the Definition \ref{mapu}.
\item [D3)] $\widehat{U}_{\phi}(w,h)\in \mathcal{U}$ as the control strategy derived from $\widehat{u}(w,h)$ through the Definition \ref{mapU}.
\end{enumerate}
\end{definition}
Notice that the existence of $\widehat{u}(w,h)$ follows from Assumptions \ref{A1}-\ref{Hyp5a} and Theorem 3.3.5 in \cite{hernandez96}, and
the fact that $\widehat{u}_{\phi}(w,h)\in \mathcal{S}_{\mathbb{V}}$, and $\widehat{U}_{\phi}(w,h)\in \mathcal{U}$ follow from Proposition 3.10 in \cite{average}.

\bigskip
\noindent
In the next assumption notice that for any $u\in \mathcal{S}_\mathbb{U}$, $G(x,u_\phi;.)$ can be seen as the stochastic kernel associated to the post-jump location of a PDMP.
This assumption is related to some geometric ergodic properties of the operator $G$ (see for example the comments on page 122 in \cite{hernandez99} or Lemma 3.3 in \cite{guo06}
for more details on this kind of assumption).

\begin{assumption}
\label{A2}
Suppose that there exist $a> 0$, $0<\kappa<1$ and for any $u\in \mathcal{S}_\mathbb{U}$ there exists a probability measure $\nu_u$, such that $\nu_u(g)<+\infty$ and
\begin{eqnarray}
\bigl| G^kh(x,u_\phi) - \nu_u(h) \bigr| \leq a \|h\|_{g} \kappa^k g(x),
\label{Hy3}
\end{eqnarray}
for all $h\in \mathbb{B}_{g}(E)$ and $k\in\NN$.
\end{assumption}
The final assumption is:

\begin{assumption}
\label{Hyp8a} There exist $\underline{\lambda} \in\mathbb{M}(E)^{+}$, $\widebar{f} \in\mathbb{M}(E)^{+}$, $K_{\lambda} \in \RR_{+}$ such that
\begin{enumerate}
\item [a)] $\lambda(y,a) \geq \underline{\lambda}(y)$ and $f(y,a)\leq \widebar{f}(y)$ for all $y\in E$ and $a\in\mathbb{U}(y)$,
\item [b)] $\ds \int_0^{t_{*}(x)} e^{ct-\int_0^t \underline{\lambda}(\phi(x,s))ds}dt \leq K_{\lambda}$, for all $x\in E$,
\item [c)] $\ds \lim_{t\rightarrow +\infty} e^{ct-\int_0^{t}\underline{\lambda}(\phi(x,s))ds} =0$, for all $x\in E$ with $t_{*}(x)=+\infty$,
\item [d)] $\ds \lim_{t\rightarrow +\infty} e^{-\int_0^t \underline{\lambda}(\phi(x,s))ds} g(\phi(x,t)) =0$, for all $x\in E$ with $t_{*}(x)=\infty$,
\item [e)] $\ds \int_0^{t_{*}(x)} e^{-\int_0^t \underline{\lambda}(\phi(x,s))ds} \widebar{f}(\phi(x,t)) dt < \infty$.
\end{enumerate}
\end{assumption}

\begin{remark} \label{vac} Notice the following consequences of Assumption \ref{Hyp8a}:
\begin{enumerate}
\item [i)] Assumption \ref{Hyp8a} c) implies that $\ds G_{\alpha}(x,\Theta;A) = \int_0^{t_{*}(x)}e^{-\alpha s - \Lambda^{\mu}(x,s)}\lambda QI_{A}(\phi(x,s),\mu(s)) ds$, and $H_{\alpha}w(x,\Theta) =0$,
for any $x\in E$ with $t_{*}(x)=+\infty$, $A\in \mathcal{B}(E)$, $\alpha\geq -c$, $\Theta=(\mu,\mu_{\partial})\in \mathbb{V}^{r}(x)$, $w\in \mathbb{M}(\partial E\times \mathbb{U})$.
\item [ii)] Assumptions \ref{Hyp8a} a) and b) imply that $\ds \mathcal{L}_{\alpha}(x,\Theta) \leq  K_{\lambda}$ for any $\alpha\geq -c$, $x\in E$, $\Theta\in \mathbb{V}^{r}(x)$.
\end{enumerate}
\end{remark}

\subsection{Properties of the $\alpha$-discount value function $\mathcal{J}_D^\alpha(\cdot)$}
\label{Value}
The next two propositions establish a connection between a general intro-differential inequality (respectively equality) related to the local characteristics of the PDMP and an inequality (respectively equality) related to the operators
$G_{\alpha}$, $L_{\alpha}$ and $H_{\alpha}$. They will be crucial for the boundedness results on $\mathcal{J}_D^\alpha(\cdot)$ to be developed in the sequel.
\begin{proposition}
\label{prop1}
Suppose that 
there exist $v\in\mathbb{M}^{ac}(E,\RR_{+})$, $\ell\in\mathbb{M}(E)^{+}$, $k\in\mathbb{M}(E)^{+}$,
$p \in\mathbb{M}(\partial E)^{+}$, $\Theta=(\mu,\mu_{\partial})\in \mathcal{S}_{\mathbb{V}^{r}}$, $d\geq 0$, and $\alpha \geq  -c$ satisfying
\begin{align}
\mathcal{X}v(\phi(x,t))-\left[\alpha+\lambda(\phi(x,t),\mu(x,t))\right] & v(\phi(x,t)) + \ell(\phi(x,t)) \nonumber \\
& +\lambda(\phi(x,t),\mu(x,t))Qk(\phi(x,t),\mu(x,t)) \leq d,
\label{eqp1}
\end{align}
for all $x\in E$, $t\in[0,t_{*}(x))$ and
\begin{eqnarray}
v(\phi(x,t_{*}(x)))\geq p(\phi(x,t_{*}(x)))+Qk(\phi(x,t_{*}(x)),\mu_{\partial}(\phi(x,t_{*}(x)))),
\label{eqp2}
\end{eqnarray}
for all $x \in E$ with $t_{*}(x)<\infty$.
\nl
Then
\begin{eqnarray}
v(x)  & \geq  & -d \mathcal{L}_{\alpha}(x,\Theta(x)) + L_{\alpha}\ell(x,\Theta(x)) + H_{\alpha}p(x,\Theta(x)) + G_{\alpha}k(x,\Theta(x)).
\label{intg}
\end{eqnarray}
\end{proposition}
\noindent {\textbf{Proof}:} Multiplying both sides of equation (\ref{eqp1}) by $\ds e^{-\alpha t - \Lambda^{\mu(x)}(x,t)}$ and integrating over $[0,s]$ for $s\in[0,t_{*}(x))$ we get that
\begin{align}
 d \int_{0}^{s} e^{-\alpha t - \Lambda^{\mu(x)}(x,t)} dt  \geq  & \: e^{-\alpha s - \Lambda^{\mu(x)}(x,s)} v(\phi(x,s))-v(x)  + \int_{0}^{s} e^{-\alpha t - \Lambda^{\mu(x)}(x,t)} \bigl[\ell(\phi(x,t)) \nonumber \\
& + \lambda(\phi(x,t),\mu(x,t))Qk(\phi(x,t),\mu(x,t))\bigr] dt
\label{Cu1a}.
\end{align}
Consider the case in which $t_{*}(x)< \infty$.
By using the fact that $v\in\mathbb{M}^{ac}(E)$, we obtain from Remark \ref{vac} $ii)$ and equation (\ref{Cu1a}) that
\begin{align}
v(x) \geq  & -d \mathcal{L}_{\alpha}(x,\Theta(x)) + L_{\alpha}\ell(x,\Theta(x))
+ e^{-\alpha t_{*}(x) -\Lambda^{\mu(x)}(x,t_{*}(x))} v(\phi(x,t_{*}(x)))\nonumber \\
& + \int_{0}^{t_{*}(x)} e^{-\alpha t - \Lambda^{\mu(x)}(x,t)} \lambda(\phi(x,t),\mu(x,t))Qk(\phi(x,t),\mu(x,t)) dt
\label{Cu1abis}.
\end{align}
However, from equation (\ref{eqp2}), it follows that
\begin{align*}
v(x) \geq  & -d \mathcal{L}_{\alpha}(x,\Theta(x)) + L_{\alpha}\ell(x,\Theta(x)) + H_{\alpha}p(x,\Theta(x)) + G_{\alpha}k(x,\Theta(x)).
\end{align*}
Now consider the case in which $t_{*}(x)=+\infty$. From equation (\ref{Cu1a}) (and recalling that $v$ is positive), we have that
\begin{align*}
 d \int_{0}^{s} e^{-\alpha t - \Lambda^{\mu(x)}(x,t)} dt  \geq   -v(x)  + \int_{0}^{s} & e^{-\alpha t - \Lambda^{\mu(x)}(x,t)} \bigl[ \ell(\phi(x,t)) \nonumber \\
 & + \lambda(\phi(x,t),\mu(x,t))Qk(\phi(x,t),\mu(x,t))\bigr] dt,
\end{align*}
and so, by taking the limit as $s$ tends to infinity in the previous equation, it yields
\begin{eqnarray*}
v(x) & \geq & -d \mathcal{L}_{\alpha}(x,\Theta(x)) + L_{\alpha}\ell(x,\Theta(x)) \nonumber\\
& & + \int_{0}^{t_{*}(x)} e^{-\alpha t - \Lambda^{\mu(x)}(x,t)} \lambda(\phi(x,t),\mu(x,t))Qk(\phi(x,t),\mu(x,t)) dt.
\end{eqnarray*}
However, by using the fact that $t_{*}(x)=+\infty$ and Remark \ref{vac} $i)$, we have  that  $H_{\alpha}p(x,\Theta(x))=0$ and
$\ds G_{\alpha}k(x,\Theta(x))=\int_{0}^{t_{*}(x)} e^{-\alpha t - \Lambda^{\mu(x)}(x,t)} \lambda(\phi(x,t),\mu(x,t))Qk(\phi(x,t),\mu(x,t)) dt$, showing the result.
\hfill$\Box$

\bigskip

\noindent
If the inequalities in (\ref{eqp1}) and (\ref{eqp2}) are replaced by equalities then the hypotheses of Proposition \ref{prop1} must be restricted to $\alpha\geq 0$ to
show that the inequality in (\ref{intg}) becomes an equality, more specifically, we have the following result:
\begin{proposition}
\label{prop2}
Suppose that 
there exist $v\in\mathbb{M}_{w}^{ac}(E,\RR_{+})$, $\ell\in\mathbb{M}(E)^{+}$, $k\in\mathbb{M}(E)^{+}$,
$p \in\mathbb{M}(\partial E)^{+}$, $\Theta=(\mu,\mu_{\partial})\in \mathcal{S}_{\mathbb{V}^{r}}$, $d\geq 0$, and $\alpha \geq  0$ satisfying
\begin{align}
\mathcal{X}v(\phi(x,t))-\left[\alpha+\lambda(\phi(x,t),\mu(x,t))\right] & v(\phi(x,t)) + \ell(\phi(x,t)) \nonumber \\
& +\lambda(\phi(x,t),\mu(x,t))Qk(\phi(x,t),\mu(x,t)) = d,
\label{eqp1=}
\end{align}
for all $x\in E$, $t\in[0,t_{*}(x))$ and
\begin{eqnarray}
v(\phi(x,t_{*}(x))) = p(\phi(x,t_{*}(x)))+Qk(\phi(x,t_{*}(x)),\mu_{\partial}(\phi(x,t_{*}(x)))),
\label{eqp2=}
\end{eqnarray}
for all $x \in E$ with $t_{*}(x)<\infty$.
\nl
Then
\begin{eqnarray}
v(x)  & =  & -d \mathcal{L}_{\alpha}(x,\Theta(x)) + L_{\alpha}\ell(x,\Theta(x)) + H_{\alpha}p(x,\Theta(x)) + G_{\alpha}k(x,\Theta(x)).
\label{intk=}
\end{eqnarray}
\end{proposition}
\noindent {\textbf{Proof}:} By following the same steps as in the first part of the proof of Proposition \ref{prop1} we have that
for all $s\in[0,t_{*}(x))$,
\begin{align}
d\int_{0}^{s} e^{-\alpha t - \Lambda^{\mu(x)}(x,t)} dt  =  & \: e^{-\alpha s - \Lambda^{\mu(x)}(x,s)} v(\phi(x,s))-v(x)  + \int_{0}^{s} e^{-\alpha t - \Lambda^{\mu(x)}(x,t)} \bigl[\ell(\phi(x,t)) \nonumber \\
& + \lambda(\phi(x,t),\mu(x,t))Qk(\phi(x,t),\mu(x,t))\bigr] dt
\label{Cu1a=}.
\end{align}
The case in which $t_{*}(x)< \infty$ can be treated in the same manner as in the proof of Proposition \ref{prop1}.
However, the case in which $t_{*}(x)=+\infty$ is different.
By using Assumption \ref{Hyp8a} d) and the fact that $0\leq v \leq \|v\|_{g} g$, we have that for any $\alpha \geq 0$,
\begin{eqnarray*}
\lim_{s\rightarrow +\infty} e^{-\alpha s-\Lambda^{\mu(x)}(x,s)} v(\phi(x,s)) \leq
\|v\|_{g} \lim_{s\rightarrow +\infty} e^{-\int_{0}^{t_{*}(x)}  \underline{\lambda}(\phi(x,t))dt} g(\phi(x,s)) = 0.
\end{eqnarray*}
Therefore, taking the limit as $s$ tends to infinity in equation (\ref{Cu1a=}), we have that
\begin{align*}
d \mathcal{L}_{\alpha}(x,\Theta(x)) =  &  -v(x)  + L_{\alpha}\ell(x,\Theta(x))  \nonumber \\
 & + \int_{0}^{s} e^{-\alpha t - \Lambda^{\mu(x)}(x,t)} \lambda(\phi(x,t),\mu(x,t))Qk(\phi(x,t),\mu(x,t)) dt,
\end{align*}
and this shows equation (\ref{intk=}) by using Remark \ref{vac} $i)$.
\hfill$\Box$

\bigskip

Applying Proposition \ref{prop1} to the inequalities (\ref{Cu1}) and (\ref{Cu2}) we obtain the following corollary:
\begin{corollary}
\label{corol1}
For any $u\in \mathcal{S}_\mathbb{U}$, $\alpha \geq -c$, and $x\in E$
\begin{eqnarray}
g(x)  & \geq  & -b \mathcal{L}_{\alpha}(x,u_{\phi}(x)) + (c+\alpha) L_{\alpha}g(x,u_{\phi}(x)) + H_{\alpha}\overline{r}(x,u_{\phi}(x)) + G_{\alpha}g(x,u_{\phi}(x)),
\label{intw}
\end{eqnarray}
and for all $\Theta \in \mathcal{S}_{\mathbb{V}^{r}}$
\begin{eqnarray}
(c+\alpha)L_{\alpha}g(x,\Theta(x)) + H_{\alpha}\overline{r}(x,\Theta(x)) + G_{\alpha}g(x,\Theta(x)) \leq bK_{\lambda}+g(x).
\label{intLHG}
\end{eqnarray}
\end{corollary}
\noindent {\textbf{Proof}:}
Clearly from Proposition 3.8 and Remark 3.11 in \cite{average}, it follows that $u_{\phi}\in \mathcal{S}_{\mathbb{V}^{r}}$.
Consequently, setting $d=b$, $v=g$, $\ell=(c+\alpha)g$, $p=\overline{r}$, $k=g$, and $\Theta=u_{\phi}$ in Proposition \ref{prop1} we get equation (\ref{intw}).
Similarly, from Remark \ref{vac} $ii)$, the inequality (\ref{intLHG}) is a straightforward consequence of the inequality (\ref{intg}).
\hfill$\Box$

\bigskip

The next theorem provides bounds in terms of $\alpha$ and $g$ for a sequence of functions defined by a general recursive equation and for the functions $Lf$, $Hr$ and $Lg$.
\begin{theorem}
\label{theo1}
Define the sequence $(q_m(x))_{m\in \NN}$ by
\begin{align}
q^\alpha_{0}(x)  & = 0, \nonumber \\
q^\alpha_{m+1}(x)  & = L_\alpha f(x,u^{m+1}_{\phi}(x)) + H_\alpha r(x,u^{m+1}_{\phi}(x))+ G_\alpha q^\alpha_{m}(x,u^{m+1}_{\phi}(x)),
\label{eqconv1}
\end{align}
where $x\in E$, $(u^{m})_{m\in \NN}\in \mathcal{S}_\mathbb{U}$ and $\alpha > 0$.
\nl
Then the following assertions hold:
\begin{enumerate}
\item [i)] for any $x\in E$, $m\in \NN$ and $\alpha \in [0, \delta)$, we have that
\begin{equation}
q^\alpha_m(x) \leq \frac{M}{c+\alpha} g(x) + \frac{Mb}{c \alpha}.
\label{Cu4}
\end{equation}
\item [ii)] for any $x\in E$, $u\in \mathcal{S}_\mathbb{U}$,
\begin{align}
0&\leq Lf(x,u_\phi(x))+Hr(x,u_\phi(x)) \leq \frac{M(1+bK_{\lambda})}{c}g(x),
\label{Cu4a}\\
0 &\leq Lg(x,u_\phi(x))\leq \frac{(1+bK_{\lambda})}{c}g(x).
\label{Cu4b}
\end{align}
\end{enumerate}
\end{theorem}
\noindent {\textbf{Proof}:}
Let us show (\ref{Cu4}) by induction. For $m=0$ it is immediate since $q_0^\alpha=0$. Suppose it holds for $m$.
Combining (\ref{eqconv1}) and  (\ref{Cu4}) we have
\begin{equation}
q^\alpha_{m+1}(x) \leq L_{\alpha}f(x,u^{m}_{\phi}(x))+H_{\alpha}r(x,u^{m}_{\phi}(x)) + \frac{M}{c+\alpha}G_{\alpha}g(x,u^{m}_{\phi}(x)) + \frac{Mb}{c \alpha} G_{\alpha}1(x,u^{m}_{\phi}(x)).
\label{eqconv2}
\end{equation}
Moreover, from equations (\ref{intw}) and (\ref{intLHG}), we obtain that
\begin{equation}
G_{\alpha}g(x,u^{m}_{\phi}(x)) \leq g(x) + b\mathcal{L}_\alpha(x,u^{m}_{\phi}(x)) - (c+\alpha)L_\alpha g(x,u^{m}_{\phi}(x)) - H_\alpha \overline{r}(x,u^{m}_{\phi}(x)).
\label{Cu5a}
\end{equation}
Replacing (\ref{Cu5a}) into (\ref{eqconv2}) and using (\ref{Cu3}) and (\ref{Cu3a}), we get
\begin{align}
q^\alpha_{m+1}(x) &\leq  L_{\alpha}(f-M g)(x,u^{m}_{\phi}(x))+H_{\alpha}(r-\frac{M}{c+\alpha}\overline{r})(x,u^{m}_{\phi}(x))+ \frac{M}{c+\alpha}g(x)\nonumber\\&+
Mb\Bigl(\frac{1}{c \alpha}G_{\alpha}1(x,u^{m}_{\phi}(x)) + \frac{1}{c+\alpha}\mathcal{L}_\alpha(x,u^{m}_{\phi}(x))\Bigr)\nonumber \\
&\leq \frac{M}{c+\alpha}g(x) +
\frac{Mb}{c\alpha}\Bigl(G_{\alpha}1(x,u^{m}_{\phi}(x)) + \alpha \mathcal{L}_\alpha(x,u^{m}_{\phi}(x))\Bigr)
\nonumber \\
&\leq \frac{M}{c+\alpha}g(x) +
\frac{Mb}{c\alpha}
\label{eqconv3}
\end{align}
since that $G_{\alpha}1(x,u^{m}_{\phi}(x)) + \alpha \mathcal{L}_\alpha(x,u^{m}_{\phi}(x)) = 1$.

\bigskip

\noindent Let us show now (\ref{Cu4a}) and (\ref{Cu4b}). For $\alpha = 0$ it follows from Remark \ref{vac} $ii)$ and equation (\ref{intw}) that
\begin{eqnarray}
g(x)+bK_{\lambda} \geq  g(x) + b\mathcal{L}(x,u_\phi(x)) \geq cLg(x,u_\phi(x)) + H \overline{r}(x,u_\phi(x))+ G g(x,u_\phi(x)),
\label{tutut}
\end{eqnarray}
showing equation (\ref{Cu4b}) since $g\geq 1$ and $\overline{r}\geq 0$.
Now, combining equations (\ref{Cu3}), (\ref{Cu3a}) and (\ref{tutut}) we get (\ref{Cu4a}), showing the last part of the result.
\hfill$\Box$

\bigskip

\noindent
Based on the previous result, we obtain the following corollary showing that the $\alpha$-discount value function
$\mathcal{J}_D^\alpha(\cdot)$ belongs to $\mathbb{B}_{g}(E)$ and providing a bound for $\mathcal{J}_D^\alpha(x)$ in terms of $\alpha$.
\begin{corollary}
\label{JD}
For any $\alpha > 0$ and $x\in E$,
\begin{equation}
\mathcal{J}_D^\alpha(x) \leq \frac{M}{c+\alpha} g(x) + \frac{Mb}{c \alpha}.\label{Jalpha}
\end{equation}
\end{corollary}
\noindent {\textbf{Proof}:}
By using Propositions 7.1 and 7.5 in \cite{average}, it can be shown that there exists $u_{\phi}^m \in \mathcal{S}_{\mathbb{V}}$ such that the sequence $\big(v^\alpha_{m}(x)\big)_{m\in \NN}$ defined by
$v^\alpha_{m+1}(x) = L_{\alpha}f(x,u_{\phi}^{m}(x))+H_{\alpha}r(x,u_{\phi}^{m}(x))+G_{\alpha}v^\alpha_m(x,u_{\phi}^{m}(x))$ and $v^\alpha_{0}(x)=0$
satisfies $v^\alpha_{m+1}\uparrow \mathcal{J}_D^\alpha(x)$ as $m \uparrow \infty$.
Therefore, considering $q^\alpha_m = v^\alpha_m$ in Theorem \ref{theo1} and taking the limit as $m \uparrow \infty$ we get (\ref{Jalpha}).
\hfill $\Box$

\bigskip

\noindent
The following technical lemma shows that $\mathcal{J}_D^\alpha(x)$ can be written as an infinite sum of iterates of the stochastic kernel $G_{\alpha}$.
Using this result, $\mathcal{J}_D^\alpha(x)$ is characterized in terms of the Markov kernel $G$ in Proposition \ref{Lem1}. This is an important property. Indeed, by using classical hypotheses on $G$
such as the geometric ergodic condition in Assumption \ref{A2}, it will be shown in Theorem \ref{theo2} that the mapping defined by
$\mathcal{J}_D^\alpha(\cdot)$-$\mathcal{J}_D^\alpha(y)$ for $y$ fixed in $E$ belongs to $\mathbb{B}_{g}(E)$.
\begin{lemma}
\label{Con1} For each $\alpha>0$ there exists $u^\alpha \in \mathcal{S}_\mathbb{U}$ such that
\begin{equation}
\mathcal{J}_D^\alpha(x) = \sum_{k=0}^\infty G_\alpha^k(L_\alpha f+H_\alpha r)(x,u^\alpha_{\phi}(x)).\label{Jalpha1}
\end{equation}
\end{lemma}
\noindent {\textbf{Proof}:} As shown in \cite[Theorem 7.5]{average}, $\mathcal{J}_{\mathcal{D}}^{\alpha}\in \mathbb{\mathbb{M}(E)}$ and
$\mathcal{J}_{\mathcal{D}}^{\alpha}(x) = \mathcal{R}_\alpha(0,\mathcal{J}_{\mathcal{D}}^{\alpha})(x)$.
Moreover, from Theorem 6.4 in \cite{average}, there exists $u^\alpha \in \mathcal{S}_\mathbb{U}$ such that the ordinary feedback measurable selector
$u^\alpha_{\phi}\in \mathcal{S}_{\mathbb{V}}$ satisfies
\begin{align}
\mathcal{J}_D^\alpha(x)& = \mathcal{R}_{\alpha}(0,\mathcal{J}_D^\alpha)(x)  = \mathcal{T}_{\alpha}(0,\mathcal{J}_D^\alpha)(x) =
L_{\alpha}f(x,u^\alpha_{\phi})(x)+H_{\alpha}r(x,u^\alpha_{\phi})+G_{\alpha}\mathcal{J}_{\mathcal{D}}^{\alpha}(x,u^\alpha_{\phi}).
\label{Iter1}
\end{align}
Iterating (\ref{Iter1}) and recalling that $\mathcal{J}_{\mathcal{D}}^{\alpha}(y)\geq 0$ for every $y$, yields for every $m\in \NN$ that,
\begin{equation}
\mathcal{J}_D^\alpha(x) = \sum_{k=0}^{m-1}G_\alpha^k(L_\alpha f+H_\alpha r)(x,u^\alpha_{\phi}(x)) + G_\alpha^m \mathcal{J}_{\mathcal{D}}^{\alpha}(x,u^\alpha_{\phi}(x)) \geq
\sum_{k=0}^{m-1}G_\alpha^k(L_\alpha f+H_\alpha r)(x,u^\alpha_{\phi}(x)).
\label{Jalpha1aa}
\end{equation}
For the control $U_{u^\alpha_\phi}\in \mathcal{U}$ (see Definition \ref{mapU}), it is easy to show that
\begin{align}
\sum_{k=0}^{m-1}G_\alpha^k(L_\alpha f+H_\alpha r)(x,u^\alpha_{\phi}(x))=
E^{U_{u^\alpha_\phi}}_{(x,0)} \Biggl[  & \int_{0}^{T_{m}} e^{-\alpha s} f\bigl(X(s),  u(N(s),Z(s),\tau(s)) \bigr) ds \nonumber \\
& + \int_{0}^{T_{m}} e^{-\alpha s} r\bigl(X(s-),u_{\partial}(N(s-),Z(s-)) \bigr) dp^{*}(s) \Biggr] ,
\label{Jalpha1ab}
\end{align}
where $U_{u^\alpha_\phi}=\big(u,u_{\partial}\big)$
From Assumption \ref{Hypjump}, $T_{m}\rightarrow \infty$, $P^{U_{u^\alpha_\phi}}$ a.s. Therefore from the monotone convergence theorem, equation (\ref{Jalpha1ab}) implies that
$\ds \sum_{k=0}^{\infty}G_\alpha^k(L_\alpha f+H_\alpha r)(x,u^\alpha_{\phi}(x))= \mathcal{D}^\alpha(U_{u^\alpha_\phi},x)$,
and from equation (\ref{Jalpha1aa})
\begin{equation}
\mathcal{J}_D^\alpha(x) \geq  \sum_{k=0}^{\infty}G_\alpha^k(L_\alpha f+H_\alpha r)(x,u^\alpha_{\phi}(x))=\mathcal{D}^\alpha(U_{u^\alpha_\phi},x).
\label{Iter2}
\end{equation}
But since $U_{u^\alpha_\phi}\in \mathcal{U}$ and $\ds \mathcal{J}_D^\alpha(x) = \inf_{U\in\mathcal{U}}\mathcal{D}^{\alpha}(U,x)$ it is clear that
$\mathcal{D}^\alpha(U_{u^\alpha_\phi},x)\geq \mathcal{J}_D^\alpha(x)$, so that (\ref{Iter2}) yields (\ref{Iter1}). \hfill $\Box$

\bigskip

\noindent The next proposition gives a characterization of $\mathcal{J}_D^\alpha(x)$ in terms of $G$.
\begin{proposition}
\label{Lem1}
For $\alpha>0$ and $u^\alpha_{\phi}$  as in Lemma \ref{Con1}, define the sequence
$\big(s^\alpha_{m}(x)\big)_{m\in \NN}$ for $x\in E$ by $s^\alpha_{0}(x)=0$ and
$s^\alpha_{m+1}(x) = L_{\alpha}f(x,u^\alpha_{\phi}(x))+H_{\alpha}r(x,u^\alpha_{\phi}(x))+G_{\alpha}s^\alpha_m(x,u^\alpha_{\phi}(x))$.
Then
\begin{equation}
\mathcal{J}_D^\alpha(x) = \lim_{m\rightarrow \infty} \sum_{k=0}^{m} G^k(L(f-\alpha s^\alpha_{m+1-k})+H r)(x,u^\alpha_{\phi}(x)).\label{Jalpha4}
\end{equation}
\end{proposition}
\noindent {\textbf{Proof}:}
By definition for all $m\in \NN$, $s^\alpha_{m}\in \mathbb{M}(E)$ and
$s^\alpha_{m+1}(x) = \sum_{k=0}^m G_\alpha^k(L_\alpha f+H_\alpha r)(x,u^\alpha_{\phi}(x))$
and clearly from Lemma \ref{Con1}, we have that $s^\alpha_{m} \uparrow \mathcal{J}_D^\alpha$ as $m \uparrow \infty$.
Applying Lemma 9.2 in \cite{average}, it can be shown that $s^{\alpha}_{m}\in \mathbb{M}^{ac}(E)$ and
for all $x\in E$, and $t\in [0,t_{*}(x))$,
\begin{align*}
s^{\alpha}_{m+1}(x)
& =  \int_0^t e^{-\alpha s-\int_{0}^{s}  \lambda(\phi(x,\theta),u^{\alpha}(\phi(x,\theta)))d\theta}\Bigl[ f(\phi(x,s),u^{\alpha}(\phi(x,s))) \nonumber \\
& \phantom{=} + \lambda(\phi(x,s),u^{\alpha}(\phi(x,s))) Qs^{\alpha}_{m}(\phi(x,s),u^{\alpha}(\phi(x,s))) \Bigr] ds \nonumber \\
& \phantom{=} + e^{-\alpha t-\int_{0}^{t}  \lambda(\phi(x,s),u^{\alpha}(\phi(x,s)))ds} s^{\alpha}_{m+1}(\phi(x,t)),
\end{align*}
implying that
\begin{eqnarray}
\mathcal{X}s^\alpha_{m+1}(x) - \left[\alpha + \lambda(x,u^\alpha(x)) \right] s^\alpha_{m+1}(x) + f(x,u^\alpha(x))
+\lambda(x,u^\alpha(x)) Qs^\alpha_{m}(x,u^\alpha(x)) = 0. \label{Jalpha5}
\end{eqnarray}
Consider the case in which $t_{*}(x)< \infty$. Since $s^{\alpha}_{m+1}\in \mathbb{M}^{ac}(E)$, this yields that
\begin{align}
s^{\alpha}_{m+1} & (x) =  L_{\alpha}f(x,u^{\alpha}_{\phi}(x)) +e^{-\alpha t_{*}(x)-\int_{0}^{t_{*}(x)}  \lambda(\phi(x,s),u^{\alpha}(\phi(x,s)))ds}s^{\alpha}_{m+1}(\phi(x,t_{*}(x))) \nonumber \\
& +\int_{0}^{t_{*}(x)} e^{-\alpha s -\int_{0}^{s} \lambda(\phi(x,\theta),u^{\alpha}(\phi(x,\theta)))d\theta}\lambda(\phi(x,s),u^{\alpha}(\phi(x,s)))Qs^{\alpha}_{m}(\phi(x,s),u^{\alpha}(\phi(x,s))) ds.
\label{eq2lem1}
\end{align}
From Assumption \ref{Hyp3a}, we have that $e^{-\int_{0}^{t_{*}(x)}  \lambda(\phi(x,s),u^{\alpha}(\phi(x,s)))ds} >0$. Therefore, combining the definition of $s^{\alpha}_{m}(x)$ and
equation (\ref{eq2lem1}), it gives
\begin{eqnarray}
s^{\alpha}_{m+1}(\phi(x,t_{*}(x))) & =  & Qs^{\alpha}_{m}(\phi(x,t_{*}(x)),u(\phi(x,t_{*}(x))))+r(\phi(x,t_{*}(x)),u(\phi(x,t_{*}(x)))).
\label{Jalpha6}
\end{eqnarray}
Using Proposition \ref{prop2} we get from (\ref{Jalpha5}), (\ref{Jalpha6}) that
\begin{align}
s^\alpha_{m+1}(x)= L(f- \alpha s^\alpha_{m+1})(x,u_\phi^\alpha(x))+Hr(x,u_\phi^\alpha(x)) + G s^\alpha_{m}(x,u_\phi^\alpha(x)). \label{Jalpha7}
\end{align}
Iterations of (\ref{Jalpha7}) over $m$ yields (\ref{Jalpha4}).
\hfill$\Box$

\bigskip

\noindent Before showing that the mapping defined by $\mathcal{J}_D^\alpha(\cdot)$-$\mathcal{J}_D^\alpha(y)$ for $y$ fixed in $E$ belongs to $\mathbb{B}_{g}(E)$, we need to prove that
the mapping $L(f- \alpha s^\alpha_{m+1})(.,u_\phi^\alpha(.))+Hr(.,u_\phi^\alpha(.))$ belongs to $\mathbb{B}_{g}(E)$.
\begin{lemma}
\label{Lem2}
Define $M^\prime = \frac{M(1+\frac{b}{c})(1+bK_{\lambda})}{c}$. 
For $\alpha>0$, $u^\alpha_{\phi}$ as in Lemma \ref{Con1}, $s^\alpha_{m}$ as in Lemma \ref{Lem1}, and $x\in E$, we have that
\begin{equation}
\bigl| L(f- \alpha s^\alpha_{m+1})(x,u_\phi^\alpha(x))+Hr(x,u_\phi^\alpha(x)) \bigr| \leq M^{\prime} g(x).
\label{norm1}
\end{equation}
\end{lemma}
\noindent {\textbf{Proof}:}
Notice that
\begin{equation}
-\alpha Ls^\alpha_{m+1}(x,u_\phi^\alpha(x))\leq L(f- \alpha s^\alpha_{m+1})(x,u_\phi^\alpha(x))+Hr(x,u_\phi^\alpha(x)) \leq Lf(x,u_\phi^\alpha(x))+Hr(x,u_\phi^\alpha(x)).
\label{norm2}
\end{equation}
Considering $q^\alpha_{m}=s^\alpha_{m}$ in Theorem \ref{theo1} and recalling that $g\geq 1$ we get from equation (\ref{Cu4}) that
\begin{equation}
s^\alpha_m(x) \leq \frac{M}{c+\alpha} g(x) + \frac{Mb}{c \alpha} \leq \frac{M(1+\frac{b}{c})}{\alpha}g(x).
\label{Cu4ab}
\end{equation}
Therefore from (\ref{Cu4ab}) we have that $\alpha s^\alpha_m \leq M(1+\frac{b}{c}) g$ and thus, from (\ref{Cu4b}),
\begin{align}
\alpha Ls^\alpha_{m+1}(x,u_\phi^\alpha(x)) 
& \leq \frac{M(1+\frac{b}{c})(1+bK_{\lambda})}{c}g(x). \label{ineq1}
\end{align}
By combining equations (\ref{Cu4a}), (\ref{norm2}) and (\ref{ineq1}) the result follows.
$\Box$

\bigskip

\noindent Finally, it is shown that $\mathcal{J}_D^\alpha(\cdot)$-$\mathcal{J}_D^\alpha(y)$ for $y$ fixed in $E$ belongs to $\mathbb{B}_{g}(E)$.
\begin{theorem}
\label{theo2}
For any $\alpha>0$ and $(x,y)\in E^{2}$
\begin{eqnarray}
\bigl| \mathcal{J}_D^\alpha(x) - \mathcal{J}_D^\alpha(y) \bigr| \leq \frac{a M^\prime}{1-\kappa}(1+g(y))g(x). \label{norm3a}
\end{eqnarray}
\end{theorem}
\noindent {\textbf{Proof}:}
From Assumption \ref{A2} and Lemma \ref{Lem2} we get that for all $x\in E$,
\begin{equation*}
\Bigl| G^k(L(f-\alpha s^\alpha_{m+1-k})+H r)(x,u^\alpha_{\phi}(x)) - \pi_{u^\alpha}\bigl(L(f-\alpha s^\alpha_{m+1-k})+H r\bigr) \Bigr| \leq a M^\prime \kappa^k g(x).
\end{equation*}
Consequently,
\begin{align*}
 \bigg| \sum_{k=0}^m G^k(L(f-\alpha s^\alpha_{m+1-k})+H r)(x,u^\alpha_{\phi}(x)) - G^k(L(f-\alpha & s^\alpha_{m+1-k})+H r)(y,u^\alpha_{\phi}(y)) \bigg| \nonumber \\
 & \leq a M^\prime (g(x)+g(y))\frac{1-\kappa^{m+1}}{1-\kappa}.
\end{align*}
Taking the limit as $m\uparrow \infty$ in the previous equation and recalling that $g\geq 1$ we get the desired result from Proposition \ref{Lem1}.
\hfill $\Box$

\subsection{Convergence and semi-continuity results}
\label{auxil}
The main goal of this sub-section is to show that there exists an ordinary feedback measurable selector for the one-stage optimization problems. First we present in the next two
results some convergence and semi-continuity properties of $G_{\alpha}$, $H_{\alpha}$, $L_{\alpha}$ and $\mathcal{L}_{\alpha}$.
\begin{proposition}
\label{lemlsc1} Consider $\alpha \in \RR_{+}$, a non increasing sequence of positive numbers $\{\alpha_k\}$ with $\alpha_k\downarrow \alpha$ and a sequence of functions $\big(h_{k}\big)_{k\in\NN}\in \mathbb{B}_{g}(E)$ such that
there exists $K_{h}$ satisfying $\bigl| h_{k}(x) \bigr| \leq K_{h} g(x)$ for all $k$ and all $x\in E$.
Set $\ds h = \liminf_{k\rightarrow \infty} h_{k}$.
For $x\in E$, consider $\Theta_n=\bigl(\mu_n,\mu_{\partial,n}\bigr)\in \mathbb{V}^{r}(x)$ and $\Theta=\bigl(\mu,\mu_{\partial}\bigr)\in \mathbb{V}^{r}(x)$ such that $\Theta_n \rightarrow \Theta$.
We have the following results:
\begin{align*}
&a) \lim_{n\rightarrow \infty}\mathcal{L}_{\alpha_n}(x,\Theta_n) =\mathcal{L}_{\alpha}(x,\Theta),& &  b) \liminf_{n\rightarrow \infty}L_{\alpha_n}f(x,\Theta_n) \geq L_{\alpha}f(x,\Theta), \\
&c) \liminf_{n\rightarrow \infty}H_{\alpha_n}r(x,\Theta_n) \geq H_{\alpha}r(x,\Theta),& &d) \liminf_{n\rightarrow \infty} G_{\alpha_n}h_{n}(x,\Theta_n) \geq G_{\alpha}h(x,\Theta).
\end{align*}
\end{proposition}
\textbf{Proof:}
The proofs of a), b), c) are the same as in Proposition 5.7 in \cite{average}. It only remains to show d).
Set $\tilde{h}_{k}= h_{k}+K_hg$, $\tilde{h}= h+K_hg$ and $\tilde{g}_k = \inf_{j\geq k}\tilde{h}_{j}$
(thus $\tilde{g}_k\uparrow \tilde{h}$ and $\tilde{g}_k\leq \tilde{h}_{n}$ for $n\geq k$).
By hypothesis, $\tilde{g}_k(y) \geq 0$ for all $y\in E$.
We have that $\tilde{g}_k$ is the limit of a nondecreasing sequence of measurable bounded functions $\tilde{g}_{k,i} \in \mathbb{B}(E)$. Set $\lambda_m(y,a) = m\wedge\lambda(y,a)$.
From Assumptions \ref{Hyp3a} and \ref{Hyp5a}, we have that for each $k$, $i$, $m$ and $y\in E$, $\lambda_m Q\tilde{g}_{k,i}(y,.)$ is continuous on $\mathbb{U}(y)$.
Assumption \ref{Hyp8a} and the fact that for each $k,i$, $\tilde{g}_{k,i}$ is bounded above by, say $M_{k,i}$, yields that
\begin{align*}
0 &\leq \int_0^{t_{*}(x)}  e^{-\int_0^t \underline{\lambda}(\phi(x,s))ds}\sup_{a\in \mathbb{U}(\phi(x,t))}(\lambda_m Q\tilde{g}_{k,i}(\phi(x,t),a))dt \leq m M_{k,i} K_{\lambda}.
\end{align*}
Since $(\lambda_m Q \tilde{g}_{k,i})(y,a)\geq 0$ and it is continuous in $a$ we have from b) that
$\ds \liminf_{n\rightarrow \infty} L_{\alpha_n}(\lambda_m Q \tilde{g}_{k,i})(x,\Theta_n) \geq L(\lambda_m Q \tilde{g}_{k,i})(x,\Theta)$, and thus, recalling
that $\tilde{g}_{k,i}\leq \tilde{h}_{n}$ for $n\geq k$ and $\lambda_m \leq \lambda$,
\begin{eqnarray*}
\liminf_{n\rightarrow \infty} L_{\alpha_n}(\lambda Q\tilde{h}_{n})(x,\Theta_n) \geq L_{\alpha}(\lambda_m Q\tilde{g}_{k,i})(x,\Theta).
\end{eqnarray*}
From the monotone convergence theorem and taking the limit over $m,i,k$ we get that
\begin{eqnarray}
\liminf_{n\rightarrow \infty} L_{\alpha_n}(\lambda Q\tilde{h}_{n})(x,\Theta_n) \geq L_{\alpha}(\lambda Q\tilde{h})(x,\Theta).
\label{eqcon3}
\end{eqnarray}
By using the same arguments as above, it can be shown that
\begin{eqnarray}
\liminf_{n\rightarrow \infty} L_{\alpha_n}(\lambda Qg)(x,\Theta_n) \geq L_{\alpha}(\lambda Qg)(x,\Theta).
\label{eqcon3b}
\end{eqnarray}
Moreover, from equation (\ref{intLHG}), we have for any $v\in \mathbb{B}_{g}(E)$ that $\ds \bigl| G_{\alpha}v(x,\widetilde{\Theta}) \bigr|\leq \|v\|_{g}(bK_{\lambda}+g(x))$
for all $x\in E$ and $\widetilde{\Theta }\in \mathbb{V}^{r}$, and hence
\begin{eqnarray*}
\liminf_{n\rightarrow \infty} L_{\alpha_n}(\lambda Q\tilde{h}_{n})(x,\Theta_n) =
\liminf_{n\rightarrow \infty} L_{\alpha_n}(\lambda Q h_{n})(x,\Theta_n) + K_h \liminf_{n\rightarrow \infty}  L_{\alpha_n}(\lambda Qg)(x,\Theta_n).
\end{eqnarray*}
Similarly
\begin{eqnarray*}
L_{\alpha}(\lambda Q\tilde{h}_{\alpha})(x,\Theta) = L_{\alpha} (\lambda Qh_{\alpha})(x,\Theta) + K_h L_{\alpha}(\lambda Qg)(x,\Theta).
\end{eqnarray*}
By combining equations 
(\ref{eqcon3}) and (\ref{eqcon3b}) we get that $\ds \liminf_{n\rightarrow \infty} L_{\alpha_n}(\lambda Q h_{n})(x,\Theta_n) \geq L_{\alpha}(\lambda Q h)(x,\Theta)$.
Using similar arguments as above and c) we can show that
\begin{eqnarray*}
\liminf_{n\rightarrow \infty} H_{\alpha_n}h_{n}(x,\Theta_n) \geq  H h(x,\Theta),
\end{eqnarray*}
completing the proof of d).
\hfill $\Box$

\bigskip

\begin{corollary}
\label{corlsc1} For $x\in E$, and $h \in \mathbb{B}_{g}(E)$, $\mathcal{L}_\alpha(x,\Theta)$ is continuous on $\mathbb{V}^{r}(x)$ and
$G_{\alpha} h(x,\Theta)$ (respectively, $L_{\alpha} f(x,\Theta)$, $H_{\alpha} r(x,\Theta)$) is lower semicontinuous on $\mathbb{V}^{r}(x)$.
\end{corollary}
\textbf{Proof:} By taking $\alpha_k=\alpha\geq 0$, $h_{k}=h$ in Proposition \ref{lemlsc1} the results follow.
\hfill $\Box$

\bigskip

\noindent The next two technical lemmas will be used to derive the main result of this sub-section, which is Theorem \ref{theo3main}.
\begin{lemma}
\label{prop3b}
Let $\alpha\geq 0$, $\rho \in \RR_{+}$, $h\in\mathbb{B}_{g}(E)$ and set $w=\mathcal{R}_{\alpha}(\rho,h)$.
Then there exists $\hat{\Theta}\in \mathcal{S}_{\mathbb{V}^{r}}$ such that
\begin{align}
w(x) & = -\rho \mathcal{L}_{\alpha}(x,\hat{\Theta}(x)) + L_{\alpha}f(x,\hat{\Theta}(x))+H_{\alpha}r(x,\hat{\Theta}(x))
+G_{\alpha} h(x,\hat{\Theta}(x)).
\label{O2prop3b}
\end{align}
Moreover, $w\in\mathbb{M}^{ac}(E)$, and satisfies for all $x\in E$ and $t\in [0,t_{*}(x))$,
\begin{align}
w(x)  & = \inf_{ \mu \in \mathcal{V}^{r}(x)} \bigg\{ \int_0^t e^{-\alpha s-\Lambda^{\mu}(x,s)}\biggl[ -\rho + f(\phi(x,s),\mu(s))
+ \lambda Qh(\phi(x,s),\mu(s))\biggr] ds \nonumber \\
& \phantom{=} + e^{-\alpha t-\Lambda^{\mu}(x,t)} w(\phi(x,t))\bigg\}
\label{O3prop3b} \\
& =  \int_0^t e^{-\alpha s-\Lambda^{\hat{\mu}(x)}(x,s)}\biggl[ -\rho + f(\phi(x,s),\hat{\mu}(x,s))
 + \lambda Qh(\phi(x,s),\hat{\mu}(x,s))\biggr] ds \nonumber \\
& \phantom{=} + e^{-\alpha t-\Lambda^{\hat{\mu}(x)}(x,t)} w(\phi(x,t)),
\label{O4prop3b}
\end{align}
where $\hat{\Theta}(x) = (\hat{\mu}(x), \hat{\mu}_{\partial}(x))$.
\end{lemma}
\noindent {\bf{Proof:}}
From Corollary \ref{corol1}, it follows that the mapping $V$ defined on $\mathcal{K}$ by
$$ V(x,\Theta)=-\rho \mathcal{L}_{\alpha}(x,\Theta) + L_{\alpha}f(x,\Theta)+H_{\alpha}r(x,\Theta)+G_{\alpha} h(x,\Theta) $$
takes values in $\RR$.
Moreover, from Assumption \ref{Mesurability} and Proposition 3.4 in \cite{average}, it follows that $V$ is measurable.
Furthermore, by using Corollary 5.8 in \cite{average} it follows that for all $x\in E$, $V(x,.)$ is lower semicontinuous on $\mathbb{V}^{r}(x)$.
Recalling that $\mathbb{V}^{r}(x)$ is a compact subset of $\mathbb{V}^{r}$ and by using Proposition D.5 in \cite{hernandez96},
we obtain that there exists $\hat{\Theta}\in \mathcal{S}_{\mathbb{V}^{r}}$ such that equation (\ref{O2prop3b}) is satisfied.
The rest of the proof is similar to the proof of Proposition 4.2 in \cite{average} and it is therefore omitted.
\hfill $\Box$

\bigskip

\begin{lemma}
\label{lem3}
Let $\alpha\geq 0$, $\rho \in \RR_{+}$ and $h\in\mathbb{B}_{g}(E)$.
Then, for all $x\in E$
\begin{eqnarray}
\mathcal{R}_{\alpha}(\rho,h)(x) \geq -(\rho+b\|h\|_{g}) K_{\lambda} - \|h\|_{g} g(x),
\label{Nw}
\end{eqnarray}
and for all $x\in E$ such that $t_{*}(x)=\infty$ and $\Theta=(\mu,\mu_{\partial})\in\mathbb{V}^{r}(x)$
\begin{align}
-\rho \mathcal{L}_{\alpha}(x, & \Theta) + L_{\alpha}f(x,\Theta)+H_{\alpha}r(x,\Theta)+G_{\alpha}h(x,\Theta) \nonumber \\
& = \lim_{t\rightarrow +\infty} \int_0^t e^{-\alpha s-\Lambda^{\mu}(x,s)}\biggl[ -\rho + f(\phi(x,s),\mu(s)) + \lambda Qh(\phi(x,s),\mu(s))\biggr] ds .
\label{lim}
\end{align}
\end{lemma}
\noindent {\bf{Proof:}} From equation (\ref{intLHG}) we have
\begin{eqnarray}
G_{\alpha}g(x,\Theta) \leq bK_{\lambda}+g(x),
\label{eqlem3}
\end{eqnarray}
for all $x\in E$ and $\Theta \in \mathbb{V}^{r}$.
Consequently, by using equation (\ref{O2prop3b}) and the fact that $f\geq 0$ and $r\geq 0$ it follows that
$\mathcal{R}_{\alpha}(\rho,h)(x) \geq  -\rho \mathcal{L}_{\alpha}(x,\hat{\Theta}(x)) +G_{\alpha} h(x,\hat{\Theta}(x)) \geq  -(\rho+b\|h\|_{g}) K_{\lambda} - \|h\|_{g} g(x)$,
showing the first part of the result.
\nl
From Assumptions \ref{Hyp8a} a), b) and e), we have that $\ds \lim_{t\rightarrow +\infty} \int_0^t e^{-\alpha s-\Lambda^{\mu}(x,s)}\bigl[ -\rho + f(\phi(x,s),\mu(s)) \bigr] ds$ exists in $\RR$, and from
equation (\ref{eqlem3}), $\ds \lim_{t\rightarrow +\infty} \int_0^t e^{-\alpha s-\Lambda^{\mu}(x,s)} \lambda Qg(\phi(x,s),\mu(s)) ds$ exists in $\RR$.
By using the fact that $h\in \mathbb{B}_{g}(E)$, it follows that the limit on the right hand side of equation (\ref{lim}) exists. Finally, from Remark \ref{vac} $i)$ we get the last part of the result.\hfill $\Box$

\bigskip

\noindent The next result shows that for any function $h\in \mathbb{B}_{g}(E)$, the one-stage optimization operators $\mathcal{R}_{\alpha}(\rho,h)(x)$ and $\mathcal{T}_{\alpha}(\rho,h)(x)$ are equal and
that there exists an ordinary feedback measurable selector for the one-stage optimization problems associated to these operators.
\begin{theorem}
\label{theo3main}
Let $\alpha\geq 0$, $\rho \in \RR_{+}$, $h\in\mathbb{B}_{g}(E)$ and set $w=\mathcal{R}_{\alpha}(\rho,h)$.
Then $w\in\mathbb{M}^{ac}(E)$ and the ordinary feedback measurable selector $\widehat{u}_{\phi}(w,h) \in \mathcal{S}_{\mathbb{V}}$ (see item D2) of Definition \ref{defmes}) satisfies
the following one-stage optimization problems:
\begin{eqnarray}
\mathcal{R}_{\alpha}(\rho,h)(x) & = & \mathcal{T}_{\alpha}(\rho,h)(x)\nonumber \\
& = & -\rho \mathcal{L}_{\alpha}(x,\widehat{u}_{\phi}(w,h)(x)) +
L_{\alpha}f(x,\widehat{u}_{\phi}(w,h)(x))+H_{\alpha}r(x,\widehat{u}_{\phi}(w,h)(x))\nonumber \\
& & +G_{\alpha}h(x,\widehat{u}_{\phi}(w,h)(x)).
 \label{minpro}
\end{eqnarray}
\end{theorem}
\noindent \textbf{Proof:}
According to Lemma \ref{prop3b}, there exists $\hat{\Theta}\in \mathcal{S}_{\mathbb{V}^{r}}$ such that for all $x\in E$ and $t\in [0,t_{*}(x))$ we have
\begin{eqnarray}
e^{-\alpha t-\Lambda^{\hat{\mu}(x)}(x,t)}w(\phi(x,t))- w(x) & = & \int_0^t e^{-\alpha s -\Lambda^{\hat{\mu}(x)}(x,s)}\biggl[ \rho - f(\phi(x,s),\hat{\mu}(x,s)) \nonumber \\
& & - \lambda Qh(\phi(x,s),\hat{\mu}(x,s)) \biggr] ds,
\label{eq2theo3a}
\end{eqnarray}
where $\hat{\Theta}(x)=(\hat{\mu}(x), \hat{\mu}_{\partial}(x))$.
Since $w\in \mathbb{M}^{ac}(E)$, we obtain from equation (\ref{eq2theo3a}) that
\begin{align*}
\mathcal{X}w(\phi(x,t))-[\alpha+\lambda(\phi(x,t),\hat{\mu}(x,t))] w(\phi(x,t))= -f(\phi(x,t),\hat{\mu}(x,t))
-\lambda Qh(\phi(x,t),\hat{\mu}(x,t))+\rho,
\end{align*}
$\eta-a.s.$ on $[0,t_{*}(x))$, implying that
\begin{align*}
-\mathcal{X} & w(\phi(x,t))+ \alpha w(\phi(x,t)) \nonumber \\
& \geq  \inf_{\mu \in \mathcal{P} \bigl(\mathbb{U}(\phi(x,t))\bigr)}  \Bigl\{f(\phi(x,t),\mu)-\lambda(\phi(x,t),\mu)w(\phi(x,t)) + \lambda Qh(\phi(x,t),\mu) \Bigr\} - \rho.
\end{align*}
However, notice that
\begin{align*}
\inf_{\mu \in \mathcal{P} \bigl( \mathbb{U}(\phi(x,t))\bigr)}  \Bigl\{f( & \phi(x,t),\mu)-\lambda(\phi(x,t),\mu)w(\phi(x,t))+ \lambda Qh(\phi(x,t),\mu) \Bigr\} - \rho \\
& =  \inf_{a\in \mathbb{U}(\phi(x,t))} \Bigl\{f(\phi(x,t),a)- \lambda(\phi(x,t),a) \bigl[ w(\phi(x,t)) - Qh(\phi(x,t),a) \bigr] \Bigr\} - \rho.
\end{align*}
Consequently, by considering the measurable selector $\widebar{u}\in \mathcal{S}_{\mathbb{U}}$
given by $\widebar{u} = \widehat{u}(w,h)$ (see Definition \ref{defmes}, D1)), we have that
\begin{align*}
-\mathcal{X} & w(\phi(x,t))+ \alpha w(x) \nonumber \\
& \geq - \rho + f(\phi(x,t),\widebar{u}(\phi(x,t)))-\lambda(\phi(x,t),\widebar{u}(\phi(x,t))) \bigl[
w(\phi(x,t))-Qh(\phi(x,t),\widebar{u}(\phi(x,t))) \bigr],
\end{align*}
$\eta-a.s.$ on $[0,t_{*}(x))$ implying that
\begin{eqnarray*}
-\mathcal{X}w(\phi(x,t))+\alpha w(\phi(x,t)) & = & - \rho + f(\phi(x,t),\widebar{u}(\phi(x,t)))\nonumber \\
& & -\lambda(\phi(x,t),\widebar{u}(\phi(x,t))) \bigl[w(\phi(x,t))-Qh(\phi(x,t),\widebar{u}(\phi(x,t))) \bigr],
\end{eqnarray*}
$\eta-a.s.$ on $[0,t_{*}(x))$, otherwise this would lead to a contradiction with equation (\ref{O3prop3b}).
Consequently, for all $t\in [0,t_{*}(x))$ it follows that
\begin{align}
w(x) = e^{-(\alpha t+\overline{\Lambda}(x,t))} & w(\phi(x,t)) +
\int_0^t e^{-(\alpha s+\overline{\Lambda}(x,s))} \Bigl[ f(\phi(x,s),\widebar{u}(\phi(x,s))) \nonumber \\
& +\lambda(\phi(x,s),\widebar{u}(\phi(x,s))) Qh(\phi(x,s),\widebar{u}(\phi(x,s))) - \rho \Bigr] ds,
\label{espoir}
\end{align}
where we set $\ds \overline{\Lambda}(x,t)=\int_{0}^{t}\lambda(\phi(x,s),\widebar{u}(\phi(x,s))) ds$.
\nl
First consider the case in which $t_{*}(x)<\infty$. We obtain, by taking the limit as $t$ tends to $t_{*}(x)$ in the previous equation,
that the ordinary feedback measurable selector $\widehat{u}_{\phi}(w,h) \in \mathcal{S}_{\mathbb{V}}$ (see item D2) of Definition \ref{defmes}) satisfies:
\begin{align}
w(x) = & e^{-(\alpha t_{*}(x)+\overline{\Lambda}(x,t_{*}(x)))}  w(\phi(x,t_{*}(x)))
-\rho \mathcal{L}_{\alpha}(x,\widehat{u}_{\phi}(w,h)(x)) + L_{\alpha}f(x,\widehat{u}_{\phi}(w,h)(x)) \nonumber \\
& +\int_0^{t_{*}(x)} e^{-(\alpha s+\overline{\Lambda}(x,s))} \lambda(\phi(x,s),\widebar{u}(\phi(x,s))) Qh(\phi(x,s),\widebar{u}(\phi(x,s)))ds.
\label{infmin1bis}
\end{align}
Define the control $\Theta(x)$ by $(\hat{\mu}(x), \mu)$ for $\mu \in \mathcal{P}\bigl(\mathbb{U}(\phi(x,t_{*}(x)))\bigr)$.
Therefore, we have that
\begin{align}
w(x) \leq &-\rho \mathcal{L}_{\alpha}(x,\hat{\Theta}(x))+ L_{\alpha}f(x,\hat{\Theta}(x))
+\int_{0}^{t_{*}(x)} e^{-\alpha s -\Lambda^{\hat{\mu}(x)}(x,s)} \lambda Qh(\phi(x,s),\hat{\mu}(x,s)) ds \nonumber \\
& +e^{-\alpha t_{*}(x)-\Lambda^{\hat{\mu}(x)}(x,t_{*}(x))} \bigl[Qh(\phi(x,t_{*}(x)),\mu) + r(\phi(x,t_{*}(x)),\mu) \bigr].
\label{eq1coro2a}
\end{align}
From equation (\ref{O4prop3b}), we have that
\begin{align*}
w(x) = & \int_0^t e^{-\alpha s-\Lambda^{\hat{\mu}(x)}(x,s)}\biggl[ -\rho + f(\phi(x,s),\hat{\mu}(x,s))  + \lambda Qh(\phi(x,s),\hat{\mu}(x,s)) \biggr] ds \nonumber \\
& + e^{-\alpha t
-\Lambda^{\hat{\mu}(x)}(x,t)}w(\phi(x,t)).
\end{align*}
Since $w\in \mathbb{M}^{ac}(E)$, this yields that
\begin{align}
w(x)  = &-\rho \mathcal{L}_{\alpha}(x,\hat{\Theta}(x)) + L_{\alpha}f(x,\hat{\Theta}(x))
+\int_{0}^{t_{*}(x)} e^{-\alpha s -\Lambda^{\hat{\mu}(x)}(x,s)} \lambda Qh(\phi(x,s),\hat{\mu}(x,s)) ds \nonumber \\
& +e^{-\alpha t_{*}(x)-\Lambda^{\hat{\mu}(x)}(x,t_{*}(x))}w(\phi(x,t_{*}(x))).
\label{eq2coro2a}
\end{align}
From Assumption \ref{Hyp3a}, we have that $e^{-\Lambda^{\hat{\mu}(x)}(x,t_{*}(x))} >0$. Therefore, combining
equations (\ref{eq1coro2a}) and (\ref{eq2coro2a}), it gives that for all $x\in E$ and $\mu \in
\mathcal{P}\bigl(\mathbb{U}(\phi(x,t_{*}(x)))\bigr)$,
\begin{eqnarray*}
w(\phi(x,t_{*}(x))) & \leq  & Qh(\phi(x,t_{*}(x)),\mu)+r(\phi(x,t_{*}(x)),\mu).
\end{eqnarray*}
Clearly, by using equation (\ref{O2prop3b}), it can be claimed that the previous inequality becomes an equality for
$\mu=\hat{\mu}_{\partial}(x)$, implying that
\begin{eqnarray*}
w(\phi(x,t_{*}(x))) & = & \inf_{\mu\in \mathcal{P}(\mathbb{U}(\phi(x,t_{*}(x))))}\{r(\phi(x,t_{*}(x)),\mu)+Qh(\phi(x,t_{*}(x)),\mu)\} \\
& = & \inf_{a\in \mathbb{U}(\phi(x,t_{*}(x)))}\{  r(\phi(x,t_{*}(x)),a)+Qh(\phi(x,t_{*}(x)),a)\}.
\end{eqnarray*}
Consequently, we have that
\begin{align}
w(\phi(x,t_{*}(x))) & = r(\phi(x,t_{*}(x)),\widebar{u}(\phi(x,t_{*}(x))))+Qh(\phi(x,t_{*}(x)),\widebar{u}(\phi(x,t_{*}(x)))).
\label{infmin1}
\end{align}
Combining equations (\ref{infmin1bis}) and (\ref{infmin1}), it follows that
\begin{align*}
w(x) =  -\rho \mathcal{L}_{\alpha}(x,\widehat{u}_{\phi}(w,h)(x)) +L_{\alpha}f(x,\widehat{u}_{\phi}(w,h)(x))+H_{\alpha}r(x,\widehat{u}_{\phi}(w,h)(x)) +G_{\alpha}h(x,\widehat{u}_{\phi}(w,h)(x)).
\end{align*}

\bigskip

\noindent
Consider now the case in which $t_{*}(x)=\infty$.
By using equation (\ref{espoir}) and (\ref{Nw}) we obtain that
\begin{align}
w(x) \geq -e^{-(\alpha t+\overline{\Lambda}(x,t))} & \bigl[(\rho+b\|h\|_{g}) K_{\lambda} + \|h\|_{g} g(\phi(x,t))\bigr]
+\int_0^t e^{-(\alpha s+\overline{\Lambda}(x,s))} \Bigl[ f(\phi(x,s),\widebar{u}(\phi(x,s))) \nonumber \\
& +\lambda(\phi(x,s),\widebar{u}(\phi(x,s))) Qh(\phi(x,s),\widebar{u}(\phi(x,s))) - \rho \Bigr] ds.
\label{isa1}
\end{align}
However, from Assumptions \ref{Hyp8a} a) and d) we obtain that
\begin{eqnarray}
\lim_{t\rightarrow +\infty}e^{-(\alpha t+\overline{\Lambda}(x,t))} \bigl[(\rho+b\|h\|_{g}) K_{\lambda} + \|h\|_{g} g(\phi(x,t))\bigr] = 0.
\label{isa2}
\end{eqnarray}
Consequently, combining equations (\ref{lim}), (\ref{isa1}) and (\ref{isa2}), the ordinary feedback measurable selector $\widehat{u}_{\phi}(w,h) \in \mathcal{S}_{\mathbb{V}}$ satisfies:
\begin{align*}
w(x) \geq & -\rho +\mathcal{L}_{\alpha}(x,\widehat{u}_{\phi}(w,h)(x)) +L_{\alpha}f(x,\widehat{u}_{\phi}(w,h)(x))+H_{\alpha}r(x,\widehat{u}_{\phi}(w,h)(x)) \\
 & +G_{\alpha}h(x,\widehat{u}_{\phi}(w,h)(x)).
\end{align*}
By using equation (\ref{O3prop3b}) it follows that the inequality in the previous equation is in fact an equality.

\bigskip

In conclusion, since $\mathbb{V}(x) \subset \mathbb{V}^{r}(x)$ it follows that $\mathcal{R}_{\alpha}(\rho,h)(x)\leq
\mathcal{T}_{\alpha}(\rho,h)(x)$. However, we have shown that $\widehat{u}_{\phi}(w,h) \in \mathcal{S}_{\mathbb{V}}$ satisfies
\begin{eqnarray*}
\mathcal{R}_{\alpha}(\rho,h)(x)  & =  & -\rho \mathcal{L}_{\alpha}(x,\widehat{u}_{\phi}(w,h)(x)) +
L_{\alpha}f(x,\widehat{u}_{\phi}(w,h)(x))+H_{\alpha}r(x,\widehat{u}_{\phi}(w,h)(x)) \\
& & +G_{\alpha}h(x,\widehat{u}_{\phi}(w,h)(x)),
\end{eqnarray*}
which is the desired result.
\hfill $\Box$

\bigskip

\section{Main results}
\label{main3}
\noindent
It has been shown in a previous work of the authors (see Theorem 6.2 in \cite{average}) that if there exists $(\rho,h) \in \RR_{+}\times \mathbb{M}(E)$ with $h$ bounded from below satisfying
the discrete-time optimality equation $\mathcal{T}(\rho,h)(x) = h(x)$, and the technical condition $\ds \limsup_{t\rightarrow +\infty} \frac{1}{t} \limsup_{m\rightarrow +\infty } E^{U}_{(x,0)} \Bigl[  h\bigl(X(t\wedge T_{m}) \bigr)\Bigr] =0$,
for all $U\in \mathcal{U}$, then there exists an ordinary feedback optimal control strategy $\widehat{U}$ for the long run average-cost problem and morevoer $\rho  = \mathcal{J}_{\mathcal{A}}(x) = \mathcal{A}(\widehat{U},x)$.
However, it is hard to obtain a solution for the discrete-time optimality equation, $\mathcal{T}(\rho,h)(x) = h(x)$. A classical method to deal with this difficulty is to follow the so-called vanishing discount approach
in order to show that there exists $(\rho,h) \in \RR\times \mathbb{M}(E)$ with $h$ bounded from below satisfying an optimality inequality of the kind $h\geq \mathcal{T}(\rho,h)$.
By using the fact that $h$ is bounded from below, the previous inequality leads to the existence of an optimal control.
In this context, a classical hypothesis (see for example Assumption 5.4.1 in \cite[page 86]{hernandez96}) is to assume that the difference of the $\alpha$-discount value functions
$\mathcal{J}_D^\alpha(\cdot)-\mathcal{J}_D^\alpha(x_{0})$ is bounded from below.
This approach has been developed in \cite[Theorem 8.5]{average} to ensure the existence of an optimal ordinary feedback control strategy.
\nl
As shown in Theorem \ref{theo2} of sub-section \ref{Value}, the hypotheses made in sub-section \ref{AssDef} yields that the difference of the $\alpha$-discount value functions $\mathcal{J}_D^\alpha(\cdot)-\mathcal{J}_D^\alpha(x_{0})$
is not necessarily bounded from below.
This result implies the existence of a pair $(\rho,h)$ satisfying $h\geq \mathcal{T}(\rho,h)$ where $\rho \in \RR_{+}$ but with $h\in \mathbb{B}_{g}(E)$.
Consequently, the result presented in \cite{average} cannot be directly used.
The idea to overcome this difficulty is to show in Proposition \ref{prop2b} that for $\widehat{u}\in \mathcal{S}_{\mathbb{U}}$,
$\ds \limsup_{t\rightarrow +\infty} \frac{1}{t} \limsup_{m\rightarrow \infty} E^{U_{\widehat{u}_{\phi}}}_{(x,0)} \Bigl[ \mathcal{T}_{\alpha}(\rho,h)\bigl(X(t\wedge T_{m}) \bigr)\Bigr] \geq 0$
in order to obtain in Theorem \ref{main} the main result of this paper, which is the existence of an optimal ordinary feedback control strategy  for the long run average-cost problem of a PDMP.

\bigskip

\noindent
First we need the following auxiliary result:
\begin{lemma}
\label{prop5}
Consider an arbitrary $u\in \mathcal{S}_{\mathbb{U}}$ and let $u_\phi$ and $U_{u_\phi}$ be as in Definitions \ref{mapu} and \ref{mapU} respectively.
For all $x\in E$ define $\widehat{g}(x) =  -b \mathcal{L}_{-c}(x,u_{\phi}(x))+G_{-c}g(x,u_{\phi}(x))$.
Then $\widehat{g}\in \mathbb{B}_{g}(E)$ and $U_{u_{\phi}}$ satisfies
\begin{align}
E^{U_{u_{\phi}}}_{(x,0)} \Bigl[  \widehat{g}\bigl(X(t\wedge T_{m}) \bigr)\Bigr] \leq e^{-ct}g(x)+\frac{b}{c}\bigl[1-e^{-ct}\bigr] + a\|\widehat{g}\|_{g} g(x) \kappa^{m} +\|\widehat{g}\|_{g}\nu_{u}(g) +b K_{\lambda}.
\label{ineqce}
\end{align}
\end{lemma}
\noindent \textbf{Proof:}
From (\ref{intw}) with $\alpha = -c$ and recalling that $\overline{r}(z)\geq 0$
we obtain that
$-b \mathcal{L}_{-c}(x,u_{\phi}(x)) +G_{-c}g(x,u_{\phi}(x)) \leq g(x)$.
Clearly, $\widehat{g}\in \mathbb{M}(E)$ is bounded from below by $-bK_{\lambda}$ from Assumption \ref{Hyp8a} b) and thus $\widehat{g}\in \mathbb{B}_{g}(E)$.
Since $\widehat{g}\in \mathbb{M}(E)$ is bounded from below, it is easy to show that
$-b E^{U_{u_{\phi}}}_{(x,0)} \Bigl[ \int_{0}^{t\wedge T_{m}} e^{cs} ds\Bigr]
+ E^{U_{u_{\phi}}}_{(x,0)} \Bigl[ e^{c(t\wedge T_{m})} \widehat{g}\bigl(X(t\wedge T_{m}) \bigr)\Bigr] \leq g(x)$,
by using the same arguments as in the proof of Proposition 4.4 in \cite{average}.
Combining Fatou's Lemma and Assumption \ref{Hypjump} we obtain that
\begin{align}
E^{U_{u_{\phi}}}_{(x,0)} \Bigl[  \widehat{g}\bigl(X(t) \bigr)\Bigr] \leq e^{-ct}g(x)+\frac{b}{c}\bigl[1-e^{-ct}\bigr].
\label{eqprop5}
\end{align}
Clearly, we have $E^{U_{u_{\phi}}}_{(x,0)} \Bigl[  \widehat{g}\bigl(X(t\wedge T_{m}) \bigr)\Bigr]
=E^{U_{u_{\phi}}}_{(x,0)} \Bigl[  I_{\{t<T_{m}\}}\widehat{g}\bigl(X(t) \bigr)\Bigr] + E^{U_{u_{\phi}}}_{(x,0)} \Bigl[  I_{\{t\geq T_{m}\}} \widehat{g}\bigl(X(T_{m}) \bigr)\Bigr]$.
Consequently, we get
$E^{U_{u_{\phi}}}_{(x,0)} \Bigl[  \widehat{g}\bigl(X(t\wedge T_{m}) \bigr)\Bigr] \leq E^{U_{u_{\phi}}}_{(x,0)} \Bigl[ \widehat{g}\bigl(X(t) \bigr)\Bigr] +G^{m}\widehat{g}(x,u_{\phi}(x))+b K_{\lambda}$
by recalling that $\widehat{g}$ is bounded from below by $-bK_{\lambda}$.
The result follows by using Assumption \ref{A2} and equation (\ref{eqprop5}).
\hfill$\Box$

\bigskip

\noindent
We have the following propositions showing that there exists $(\rho,h) \in \RR_{+}\times \mathbb{B}_{g}(E)$ such that the optimality inequality $h\geq \mathcal{T}(\rho,h)$ is satisfied:
\begin{proposition}
\label{lemLas}
Set $\rho_{\alpha} = \alpha \mathcal{J}_{\mathcal{D}}^{\alpha}(x_0)$ for a fixed state $x_0\in E$.
Then there exists a decreasing sequence of positive numbers $\alpha_k\downarrow 0$ such that $\rho_{\alpha_k}\rightarrow \rho$ and for all $x\in E$,
$\lim_{k\rightarrow \infty} \alpha_k  \mathcal{J}^{\alpha_k}(x) = \rho.$
\end{proposition}
\textbf{Proof:} From equation (\ref{Jalpha}), we obtain that there exists $\beta>0$, $C \geq 0$, such that for $\alpha \in (0,\beta]$,
$\rho_{\alpha} \leq C$. By using the lemma on page 88 in \cite{hernandez96}, the result follows.
\hfill $\Box$

\begin{proposition}
\label{lemlsc2}
Set $h_{\alpha} (\cdot)= \mathcal{J}_{\mathcal{D}}^{\alpha}(\cdot)-\mathcal{J}_{\mathcal{D}}^{\alpha}(x_0)$ for $x_0\in E$ as in Proposition \ref{lemLas} and write
$\ds h = \liminf_{k\rightarrow \infty} h_{\alpha_k}$. Then for all $x\in E$, $h \in \mathbb{B}_{g}(E)$ and $h(x)\geq\mathcal{T}(\rho,h)(x)$.
\end{proposition}
\textbf{Proof:}
From Proposition 7.1 and Theorem 7.5 in \cite{average} we have that the following equation is satisfied for each $\alpha > 0$ and $x\in E$:
\begin{align}
h_\alpha(x) & =
\mathcal{T}_\alpha(\rho_\alpha,h_\alpha)(x)\nonumber\\&
= -\rho_{\alpha}\mathcal{L}_{\alpha}(x,u^{\alpha}_{\phi}(x))+ L_{\alpha}f(x,u^{\alpha}_{\phi}(x))+H_{\alpha}r(x,u^{\alpha}_{\phi}(x))+G_{\alpha}h_{\alpha}(x,u^{\alpha}_{\phi}(x)),
\label{eqhal}
\end{align}
for $u^{\alpha}_{\phi} \in \mathcal{S}_{\mathbb{V}}$.
For $x\in E$ fixed and for all $k\in \NN$, $u^{\alpha_{k}}_{\phi}(x)\in \mathbb{V}(x)\subset \mathbb{V}^{r}(x)$ and since $\mathbb{V}^{r}(x)$ is compact we can
find a further subsequence, still written as $u^{\alpha_{k}}_{\phi}(x)$ for notational simplicity, such that $u^{\alpha_{k}}_{\phi}(x) \rightarrow \hat{\Theta} \in \mathbb{V}^{r}(x)$.
Combining equations (\ref{norm3a}), (\ref{eqhal}) and Proposition \ref{lemlsc1}
\begin{align}
\label{seleha}
h(x) & = \liminf_{k\rightarrow \infty}\Bigl\{-\rho_{\alpha_k} \mathcal{L}_{\alpha_k}(x,u^{\alpha_{k}}_{\phi}(x)) + L_{\alpha_k}f(x,u^{\alpha_{k}}_{\phi}(x))
+H_{\alpha_k}r(x,u^{\alpha_{k}}_{\phi}(x))+G_{\alpha_k}h_{\alpha_k}(x,u^{\alpha_{k}}_{\phi}(x))\Bigr\} \nonumber\\
& \geq -\rho\mathcal{L}(x,\hat{\Theta}) + L f(x,\hat{\Theta})+Hr(x,\hat{\Theta})+Gh(x,\hat{\Theta}).
\end{align}
Therefore, from Theorem \ref{theo3main}, it follows that
\begin{align*}
h(x) & \geq \mathcal{R}(\rho,h)(x) = \mathcal{T}(\rho,h)(x)
\end{align*}
showing the result.
\hfill $\Box$

\bigskip

From now on, $\rho$ and $h$ are fixed as in Propositions \ref{lemLas} and \ref{lemlsc2}, and set
$\widehat{u}=\widehat{u}(\mathcal{T}(\rho,h),h)$.
Clearly, it satisfies the following one-stage optimization problems:
\begin{eqnarray}
\mathcal{R}(\rho,h)(x) & = & \mathcal{T}(\rho,h)(x)\nonumber \\
& = & -\rho \mathcal{L}(x,\widehat{u}_{\phi}(x)) +
Lf(x,\widehat{u}_{\phi}(x))+Hr(x,\widehat{u}_{\phi}(x)) +Gh(x,\widehat{u}_{\phi}(x)).
\label{minpro2}
\end{eqnarray}

\bigskip

\noindent We need to show that
$\limsup_{t\rightarrow +\infty} \frac{1}{t} \limsup_{m\rightarrow \infty} E^{U_{\widehat{u}_{\phi}}}_{(x,0)} \Bigl[ \mathcal{T}(\rho,h)\bigl(X(t\wedge T_{m}) \bigr)\Bigr] \geq 0$.
The next proposition provides this result.

\begin{proposition}
\label{prop2b}
For all $x\in E$, $E^{U_{\widehat{u}_{\phi}}}_{(x,0)} \Bigl[ \mathcal{T}(\rho,h)\bigl(X(t\wedge T_{m}) \bigr)\Bigr]$ is well defined and satisfies
\begin{eqnarray}
\limsup_{t\rightarrow +\infty} \frac{1}{t} \limsup_{m\rightarrow \infty} E^{U_{\widehat{u}_{\phi}}}_{(x,0)} \Bigl[ \mathcal{T}(\rho,h)\bigl(X(t\wedge T_{m}) \bigr)\Bigr] \geq 0.
\label{eqlim}
\end{eqnarray}
\end{proposition}
\noindent \textbf{Proof:}
By definition, we have that $\mathcal{T}(\rho,h)(x)\geq -\rho \mathcal{L}(x,\widehat{u}_{\phi}(x))+Gh(x,\widehat{u}_{\phi}(x))$.
Therefore, using the definition of $\widehat{g}$ in Lemma \ref{prop5} with $u=\widehat{u}$ we obtain that
\begin{eqnarray}
\mathcal{T}(\rho,h)(x)  \geq -(\rho+b\|h\|_{g}) K_{\lambda} - \|h\|_{g} \widehat{g}(x).
\label{eqprop2}
\end{eqnarray}
Consequently, combining equations (\ref{ineqce}) and (\ref{eqprop2}) we get that the negative part of $ \mathcal{T}(\rho,h)\bigl(X(t\wedge T_{m}) \bigr)$ is integrable implying that
$E^{U_{\widehat{u}_{\phi}}}_{(x,0)} \Bigl[ \mathcal{T}(\rho,h)\bigl(X(t\wedge T_{m}) \bigr)\Bigr]$ is well defined, and that (\ref{eqlim}) holds,
showing the result.
\hfill$\Box$

\bigskip

\noindent The next theorem, which is the main result of this paper, shows that the ordinary feedback control $U_{\widehat{u}_{\phi}}$ is an optimal strategy for the long run average-cost problem of a PDMP.
\begin{theorem}
\label{main}
For all $x\in E$,
\begin{eqnarray*}
\rho  = \mathcal{J}_{\mathcal{A}}(x) = \mathcal{A}(U_{\widehat{u}_{\phi}},x).
\end{eqnarray*}
\end{theorem}
\textbf{Proof:}
Define
\begin{align*}
J^{U_{\widehat{u}_{\phi}}}_{m}(t,x) =  & E^{U_{\widehat{u}_{\phi}}}_{(x,0)} \biggl[ \int_{0}^{t\wedge T_{m}} \Bigl[ f\bigl(X(s),\hat{u}(X(s)) \bigr) -\rho \Bigr] ds \nonumber \\
& + \int_{0}^{t\wedge T_{m}}
r\bigl(X(s-),\hat{u}_{\partial}(X(s-)) \bigr) dp^{*}(s)
+ \mathcal{T}(\rho,h)\bigl(X(t\wedge T_{m}) \bigr)\biggr].
\end{align*}
From Proposition \ref{prop2b} we have that $E^{U_{\widehat{u}_{\phi}}}_{(x,0)} \Bigl[ \mathcal{T}(\rho,h)\bigl(X(t\wedge T_{m}) \bigr)\Bigr]$ is well defined.
Consequently, following the same arguments as in the proof of Proposition 4.4 in \cite{average}, we can show that $J^{U_{\widehat{u}_{\phi}}}_{m}(t,x)  \leq h(x)$ for all $m\in \NN$, $(t,x)\in \RR_{+}\times E$.
Therefore,
\begin{align*}
E^{U_{\widehat{u}_{\phi}}}_{(x,0)} \biggl[ & \int_{0}^{t\wedge T_{m}} \Bigl[ f\bigl(X(s),\hat{u}(X(s)) \bigr) \Bigr] ds
+ \int_{0}^{t\wedge T_{m}} r\bigl(X(s-),\hat{u}_{\partial}(X(s-)) \bigr) dp^{*}(s) \biggr] \nonumber \\
& + E^{U_{\widehat{u}_{\phi}}}_{(x,0)} \Bigl[ \mathcal{T}(\rho,h)\bigl(X(t\wedge T_{m}) \bigr)\Bigr] \leq  \rho \: t +h(x).
\end{align*}
Combining Assumption \ref{Hypjump}, the monotone convergence theorem and equation (\ref{eqlim}), it follows that
\begin{align*}
\limsup_{t\rightarrow \infty} \frac{1}{t} E^{U_{\widehat{u}_{\phi}}}_{(x,0)} \biggl[ & \int_{0}^{t} \Bigl[ f\bigl(X(s),\hat{u}(X(s)) \bigr) \Bigr] ds
+ \int_{0}^{t} r\bigl(X(s-),\hat{u}_{\partial}(X(s-)) \bigr) dp^{*}(s) \biggr] \nonumber \\
& \leq  \rho
\end{align*}
showing that $\mathcal{J}_{\mathcal{A}}(x)\leq \mathcal{A}(U_{\widehat{u}_{\phi}},x) \leq \rho$.
However, according to Theorem 1 in \cite[ chapter 5]{widder41} we have that
$\limsup_{\alpha \downarrow 0} \alpha \mathcal{J}_{\mathcal{D}}^{\alpha}(x) \leq \mathcal{J}_{\mathcal{A}}(x)$.
Consequently, from Proposition \ref{lemLas} it follows that $\rho \leq \mathcal{J}_{\mathcal{A}}(x)$, completing the proof.
\hfill $\Box$

\bibliography{Vanishing}

\end{document}